\documentclass[12pt]{amsart}
\pdfoutput=1 
\usepackage{bbm}
\usepackage[all,cmtip]{xy}
\usepackage{amsmath,amssymb,amsthm,mathabx,mathtools}
\usepackage{bbold}
\usepackage{thmtools}
\usepackage{mathrsfs}
\usepackage{geometry}
\usepackage[pagewise]{lineno}
\usepackage[backref=page, bookmarks=false, hidelinks]{hyperref}
\hypersetup{
    pdfauthor={Angel Toledo},
    pdftitle={Relative Monoidal Bondal-Orlov},
    colorlinks=true,
    citecolor=blue,
    linkcolor=black,
}
\newcommand{\dtee} {
\otimes_{X}^{\mathbb{L}}
}
\newcommand{\cc}[1][R]{\mathscr{C}h(#1)}

\newcommand{\mr}[1]{\mathrm{#1}}
\newcommand{\dteee} { \otimes^{\mathbb{L}} }
\newcommand{\ihom}[2]{\mathbb{R}\underline{\mr{Hom}}(#1,#2)}
\newcommand{\ihomc}[2]{\mathbb{R}\underline{\mr{Hom}}_{c}(#1,#2)}
\newcommand{\tofsh}{\otimes_{\mathcal{F}}}
\newcommand{\tofdgc}{\otimes_{dg}}
\newcommand{\bape}{\mathbb{A}_{pe}}
\newcommand{\bbpe}{\mathbb{B}_{pe}}

\newcommand{\bGamma}{\mathbb{\Gamma}}

\newcommand{\A}{\mathscr{A}}

\newcommand{\B}{\mathscr{B}}
\newcommand{\C}{\mathscr{C}}

\newcommand{\F}{\mathscr{F}}

\newcommand{\I}{\mathscr{I}}

\newcommand{\M}{\mathscr{M}}
\newcommand{\Ox}{\mathscr{O}}

\newcommand{\T}{\mathscr{T}}
\newcommand{\X}{\mathscr{X}}

\newcommand{\der}[1][X]{D^{b}(#1)}
\newcommand{\hzero}[1][\T]{H^{0}(#1)}

\theoremstyle{plain}
\newtheorem{thm}{Theorem}[section]
\newtheorem*{thm*}{Theorem} 

\newtheorem{prop}[thm]{Proposition}
\newtheorem{lemma}[thm]{Lemma}
\newtheorem*{lemma*}{Lemma}
\newtheorem{defn}[thm]{Definition}
\newtheorem*{defn*}{Definition}

\newtheorem*{conjecture*}{Conjecture}
\newtheorem{exmp}{Example}[section]
\newtheorem{obs}[thm]{Remark}
\newtheorem{cor}[thm]{Corollary}
\newtheorem*{cor*}{Corollary}

\declaretheoremstyle[
  headfont=\color{blue}\normalfont\bfseries,
  bodyfont=\color{black}\normalfont\itshape,
]{claim}

\declaretheoremstyle[headfont=\color{black}\normalfont,bodyfont=\color{red}\normalfont\itshape,]{proofC}

\declaretheoremstyle[headfont=\color{cyan}\normalfont, bodyfont=\color{red}\normalfont\itshape,]{Strat}

\numberwithin{equation}{section}
\providecommand{\keywords}[1]{\textbf{\textit{Index terms---}} #1}

\let\underbrace\LaTeXunderbrace

\begin{document}
\title{Relative Monoidal Bondal-Orlov}
\author[A. Sheshmani]{Artan Sheshmani}
\address{Artan Sheshmani}
\address{Beijing Institute of Mathematical Sciences and Applications, A6, Room 205
No. 544, Hefangkou Village, Huaibei Town, Huairou District, Beijing 101408,
China}
\address{Massachusetts Institute of Technology (MIT), IAiFi Institute, 182 Memorial Drive, Cambridge, MA 02139, USA}
\address{National Research University, Higher School of Economics, Russian Federation, Laboratory of Mirror Symmetry, NRU HSE, 6 Usacheva Street, Moscow, Russia, 119048}
\email{artan@mit.edu}
\author[A. Toledo]{Angel Toledo}
\address{Angel Toledo}
\address{Beijing Institute of Mathematical Sciences and Applications \\ Hefangkou Village, Huaibei Town, Huairou District, Beijing 101408, China}
\email{angeltoledo@bimsa.cn}

\date{}
\keywords{
Tensor triangulated geometry, dg-categories, derived categories, Bondal-Orlov reconstruction}
\begin{abstract}
In this article we study a relative monoidal version of the Bondal-Orlov reconstruction theorem. We establish a uniqueness result for tensor triangulated category structures $(\boxtimes,\mathbb{1})$ on the derived category $\der$ of a variety $X$ which is smooth projective and faithfully flat over a quasi-compact quasi-separated base scheme $S$ in the case where the fibers $X_{s}$ over any point $s$ of $S$ all have ample (anti-)canonical bundles. To do so we construct a stack $\bGamma$ of dg-bifunctors which parametrize the local homotopical behaviour of $\boxtimes$, and we study some of its properties around the derived categories of the fibers $X_{s}$.
\end{abstract}
\maketitle
\tableofcontents
\begingroup
\let\clearpage\relax
\section{Introduction}
The Bondal-Orlov reconstruction theorem is one of the fundamental results in the study of derived categories of coherent sheaves in algebraic geometry. This result says that when one takes a smooth projective variety $X$ over an algebraically closed field $k$ of char 0, then the bounded derived category $\der[X]$ contains enough information to reconstruct X whenever $X$ has an ample (anti) canonical bundle. \\
One of the important aspects of this result is that it establishes a clear link between the geometry, and in particular the birational geometry, of a variety and the information reflected in the whole derived category $\der[X]$ seen as an invariant of $X$ in itself. \\
In \cite{toledo2024tensor}, the second author established a monoidal version of this theorem (see Theorem \ref{thm:mbo}). To motivate this result we should think about Balmer's construction of the spectrum of a tensor triangulated category and his subsequent reconstruction theorem.\\
Concretely, Balmer describes in \cite{balmer2005spectrum} a construction which inputs a tensor triangulated category, which consists of a triangulated category with a compatible monoidal structure (Definition \ref{def:ttc}), and outputs a locally ringed space, now known as the Balmer spectrum of the tensor triangulated category. He also showed that when this tensor triangulated category is equivalent to the derived category of perfect complexes over a topologically noetherian scheme $X$ equipped with the usual derived tensor product, then the locally ringed space produced is isomorphic to $X$. \\
By this result one can wonder whether, under the hypothesis of the Bondal-Orlov theorem, one has a uniqueness of tensor triangulated category structures which by Balmer's construction allows us to deduce the conclusion of the theorem of Bondal and Orlov. \\
The answer to this question obtained in \cite{toledo2024tensor} is in the positive under some reasonably mild conditions. \\
\\
In this work we are interested in producing a similar result which now considers tensor triangulated category structures on the derived category of a smooth projective variety $X\to S$ faithfully flat over a quasi-compact and quasi-separated base scheme $S$. \\
A version of this result in the classical sense was proven by Calabrese in \cite{calabrese2013moduli}. To be precise, the full generality of Calabrese's theorem is, in his language, the following:
\begin{thm}[\cite{calabrese2018relative}, Theorem 6.2]
Let $S$ be a noetherian Artin stack with affine diagonal and let $X,Y\to S$ be flat, proper and relative algebraic spaces. Assume also that for all $s\in S$ the fibers $X_{s}, Y_{s}$ are projective, connected and Gorenstein and that $X_{s}$ has either ample or anti-ample canonical bundle. Then $X\cong Y$ as $S$-stacks if and only if there exists an $S$-linear Fourier-Mukai equivalence $\der\simeq \der[Y]$.
\end{thm}
Calabrese's approach consists of constructing two stacks $Pt_{X}$ and $BO_{X}$, of what he calls respectively lazily point-like objects and Bondal-Orlov points. Using the first such space he establishes an isomorphism between $X$ and $Pt_{X}/BG_{m}$ as stacks and then by reproducing the original argument in the proof of Bondal and Orlov at each fiber. He is then able to extend this local situation to a global result. \\
Our approach is roughly similar, although Calabrese's approach is in spirit a comparison via the well known Gabriel-Rosenberg reconstruction theorem (see \cite{calabrese2013moduli,brandenburg2018rosenberg}), while in our situation we are forced to make some higher categorical considerations. Let us phrase our main result and main corollary of this work.
\begin{thm}\label{thm:INTROTEMPMain2}
    Let $X\to S$ be a smooth projective variety faithfully flat over a quasi-compact quasi-separated scheme $S$ and suppose that there exists $s\in S$ such that $X_{s}$ has either an ample canonical bundle or ample anti canonical bundle. Let $(\boxtimes,\mathbb{1})$ be an $S$-linear tensor triangulated category structure on $Perf(X)$ (Definition \ref{def:ttc}) which is geometrically reasonable (Definition \ref{def:georeasonable}) and reasonable with respect to the adjunction pair $I_{s}^{\ast}\dashv I_{s}\,_{\ast}$ (Definition \ref{def:reasonable}). Then there exists an open affine subset $U\subseteq S$ such that $\boxtimes_{U}$ coincides on objects with $\dteee$.
\end{thm}
As a corollary, we obtain
\begin{cor}\label{cor:INTROTEMPMainCor}
Suppose $X\to S$ is as above, and assume that $Spec(\boxtimes)$ is a smooth projective variety, faithfully flat over $S$ such that $\der[X]\simeq \der[Spec(\boxtimes)]$ then $X\cong Spec(\boxtimes)$.
\end{cor}
We proceed as follows, in Section \ref{sec:prelim} we will give a brief review of the basic terminology and background results that we will be using. Concretely we will review the basic theory of tensor triangulated geometry due to Balmer, and we briefly recall the results in \cite{toledo2024tensor} where the second author proves Theorem \ref{thm:INTROTEMPMain2} for the case $S=\mr{Spec k}$. \\
We then move to investigate the relative situation which concerns us. In Section \ref{sec:xtofibers} we give some bounding conditions for which a tensor structure $\boxtimes$ in consideration can be expected to behave well with the underlying topological structure of $X\to S$. Under these conditions we show, by use of the main theorem of \cite{toledo2024tensor} (Theorem \ref{thm:mbo}) that one can restrict to a tensor structure at the level of fibers $X_{s}$ which will be then equivalent to the usual derived tensor structure $\dteee_{X_{s}}$ (Theorem \ref{thm:mainthm1}). \\
Next, we provide a quick sketch of the Morita theory of dg-categories as conceived by Toën in \cite{toen2007homotopy} and which when combined with the theory of dg-enhancements of derived categories and triangulated functors can be used to provide a framework at the level of dg-categories which enhances tensor triangulated category structures. We focus on the conditions relevant to our contexts and reflect on the concrete structures we are interested in. Concretely we are less interested on the associativity condition and coherence data of tensor triangulated structures, and more focused on just the data provided by the underlying bifunctor $\boxtimes$ ( See Definition \ref{def:pseudodgmagmoidal} ). \\
The content of Section \ref{sec:descent} is dedicated to recalling a descent construction for bounded chain complexes due to Hirschowitz and Simpson (\cite{hirschowitz1998descente}) and which we will heavily use as we consider an\ $\infty$-stack $\bape$ of dg-categories of perfect complexes over $S$ (Definition \ref{def:bape}) and of a particular kind of dg-functors between such dg-categories (Definition \ref{thm:descentformodcats}). \\
The purpose of this construction is two-fold. On the one hand the descent result of Hirschowitz-Simpson allows us to justify working always over some affine base scheme $S$ through the rest of the paper and we can exploit the rich categorical structure of the ($\infty$-)category of dg-categories to freely construct objects unavailable to us in the triangulated case. \\

The main part of this work is in Section \ref{sec:inftystack}, here we develop some fundamental ideas of the stack of dg-categories of a space $X$  (Definition \ref{def:bape}) and reason in an algebraic-geometric direction about its properties. Similarly we investigate the sheaf $\bGamma$ (Corollary \ref{prop:gammaisheaf}) which will parametrize the behaviour of $\boxtimes$ over each subset $X\times_{S} U$ for $U\subseteq S$.  \\
We show that the stack $\bape$ and the sheaf $\bGamma$ behave well with an abstract stalk construction (Definition \ref{defn:inftystalk}), and we see how this object relates to the behaviour of the dg-enhancement of the restriction of $\boxtimes$ to the fiber $X_{s}$ (Lemma \ref{thm:fiberisq}, Corollary \ref{cor:bgammasisq}). Our key lemma is a classical local-to-global result (Lemma \ref{lemma:thomasonsupp}) which is an extrapolation of a lemma due to Thomason in \cite{thomason1997classification}.  \\
With all of these results we can finally prove our second main Theorem \ref{thm:INTROTEMPMain2} and main Corollary \ref{cor:INTROTEMPMainCor}. \\
Our last Section \ref{sec:examples} we dedicate to briefly sketch some situations in which the general constructions of this article can be applied to, and potential future directions of research coming from algebraic geometry and mirror symmetry. 
\section{Preliminaries on tensor triangulated categories}\label{sec:prelim}
This section is devoted to provide preliminaries of the general theory of tensor triangulated categories and Balmer's spectra. 
\begin{defn}\label{def:ttc}
    Let $\T$ be a triangulated category, a tensor triangulated category structure (TTC for short) on $\T$ is the data of a symmetric monoidal category structure $(\boxtimes,\mathbb{1})$ such that $\mathbb{1}$ denotes the unit of the monoidal structure, the bifunctor
\[ \boxtimes:\T\times \T\to \T\]
is exact on both entries, and there exist natural isomorphisms
\[ X[1]\boxtimes Y\to (X\boxtimes Y)[1] \]
and
\[ X\boxtimes Y[1] \to (X\boxtimes Y)[1] \]
such that the following diagrams commute
\begin{align}\label{diagrams:ttccoherence1}
\begin{minipage}{0.25\textwidth}
\xymatrix{
X[1]\boxtimes U\ar[r] \ar[dr] &  X[1]  \\ 
& (X\boxtimes U)[1] \ar[u]
}
\end{minipage}
\vspace{0.1in}
\begin{minipage}{0.25\textwidth}
    \xymatrix{ U\boxtimes Y[1] \ar[r] \ar[dr] & Y[1] \\ & (U\boxtimes Y)[1] \ar[u]}
\end{minipage}
\end{align}
And the following diagram is anticommutative
\begin{align}\label{diagrams:ttccoherence2}
    \xymatrix{ X[1]\boxtimes Y[1] \ar[r] \ar[d]  & (X\boxtimes Y[1])[1] \ar[d] \\ (X[1]\boxtimes Y)[1] \ar[r] & (X\boxtimes Y)[2]}
\end{align}
This for every pair of objects $X,Y\in \T$.
\end{defn}
In \cite{balmer2005spectrum}, Balmer gave a general construction which inputs a TTC $(\boxtimes,\mathbb{1})$ on a triangulated category $\T$ and outputs a locally ringed space $\mr{Spec}(\boxtimes, \mathbb{1})$. We quickly review this construction here. \\
Balmer's construction is in spirit very similar to the construction of an affine scheme as given by the Zariski spectrum of a commutative ring. We need to introduce the following string of definitions:
\begin{defn}\label{def:thicksubt}
    Let $\I\subseteq\T$ be a full triangulated subcategory of $\T$, we say that $\I$ is thick if it is closed under direct summands.
\end{defn}
\begin{defn}\label{def:ttcideal}
Let $\T$ be a triangulated category equiped with a TTC $(\boxtimes, \mathbb{1})$. We say that a thick subcategory $\I\subseteq \T$ is a $\boxtimes$-ideal if $\T\boxtimes \I\subseteq \I$. \\
We say that such a $\boxtimes$-ideal is prime if $X\boxtimes Y\in\I$ implies $X\in \I$ or $Y\in \I$.
\end{defn}
As in affine algebraic geometry, we can define the spectrum of a tensor triangulated category using these definitions.
\begin{defn}\label{def:bspec}
    Let $(\boxtimes,\mathbb{1})$ be a TTC on a triangulated category $\T$, the set of all prime $\boxtimes$-ideals will be denoted by $\mr{Spec}(\boxtimes)$, alternatively $\mr{Spec}(\T,\boxtimes)$ depending on the given context.
 \end{defn}
 Importantly, whenever the triangualted category $\T$ is non-zero we have that $\mr{Spec}(\boxtimes)\not\emptyset$ for any tensor triangulated cateogry structure $\boxtimes$ we can put on $\T$ (see \cite[Proposition 2.3]{balmer2005spectrum}). \\
 To this set we can put a topology,
 \begin{defn}\label{def:ztopology}
Let $\T$ be a triangulated category with a TTC structure $(\boxtimes,\mathbb{1})$, the support of an object $A\in \T$ with respect to this TTC and denoted by $\mr{supp}(A)$ is the set $\{\mathfrak{p}\in \mr{Spec}(\boxtimes)\mid A\not\in \mathfrak{p}\}$.
\end{defn}
\begin{lemma}\cite[Lemma 2.6]{balmer2005spectrum}\label{lemma:zbase}
The set of the form $Z(S):=\bigcap_{A\in S}\mr{supp}(A)$, for a family of objects $S\subset \T$ forms a basis for a topology on $\mr{Spec}(\boxtimes)$.
\end{lemma}
 This topological space is always sober. For any subset $Y\subseteq \mr{Spec}(\boxtimes)$ we can define $\I_{Y}=\{A\in \T\mid \mr{supp}(A)\in Y\}$ and with this define a structure sheaf on $\mr{Spec}(\boxtimes)$ as the sheaffifcation of the endomorphism ring of the image of the unit $\mathbb{1}$ in the Verdier quotient $\T/\I_{Y}$. \\
 It is possible to show that the Balmer spectrum construction is functorial from a category of tensor triangulated categories and $\boxtimes$-exact functors (i.e. those exact functors which are also monoidal with respect to $\boxtimes$) and the construction lands in the category of locally ringed spaces. \\
The relevance of Balmer's construction in algebraic geometry is that it provides a reconstruction theorem for a large class of schemes. Namely he showed:
\begin{thm}(\cite{balmer2005spectrum})
Let $X$ be a topologically noetherian scheme and consider the triangulated category $Perf(X)$ equipped with the derived tensor product tensor triangulated category structure $(\dtee, \Ox_{X})$. Then $Spec(\dtee)\cong X$.
\end{thm}
In contrast with this result there exists a classical and well-known theory of Fourier-Mukai partners, varieties with triangulated equivalent derived categories but which might in principle not be isomorphic to each other. This kind of phenomenon was first encountered by Mukai when studying abelian varieties and it sparked a large body of research. Specifically our interest lies in the reconstruction theorem of Bondal and Orlov.
\begin{thm}\label{thm:bor}
Let $X$ be a smooth projective variety with either ample canonical or anti-canonical bundle. Suppose $\der[X]\simeq \der[Y]$ where $Y$ is another smooth projective variety. Then $X\cong Y$.
\end{thm}
In \cite{toledo2024tensor} the second author proves a monoidal version of this reconstruction result for varieties over $\mathbb{C}$, combining them with the Balmer reconstruction philosophy under the specific conditions in the hypothesis of the Bondal-Orlov reconstruction theorem.
\begin{thm}\label{thm:mbo}
    Let $X$ be a smooth projective variety with ample (anti)canonical bundle. Suppose $D^{b}(X)$ is equipped with a tensor triangulated category structure $(\boxtimes, \Ox_{X})$. Then $\boxtimes$ coincides on objects with $\dtee$.
\end{thm}
As a corollary we can obtain the following:
\begin{cor}
    Let $X$ and $\boxtimes$ be as above. Suppose $Spec(\boxtimes)$ is a smooth projective variety. Then $\der[X]\simeq \der[Spec(\boxtimes)]$ if and only if $X\cong Spec(\boxtimes)$. 
\end{cor}
We want to remark that this monoidal version is not superfluous in the presence of both reconstruction theorems. This is because as it is shown in \cite{liu2013recovering} by Liu and Sierra there are varieties $X$ such that one can equip $\der[X]$ with a tensor structure which not only is not equivalent to that of $\dtee$ but which under Balmer's reconstruction theorem produces a different locally ringed space. This phenomenon occurs already with projective spaces $\mathbb{P}^{n}$ for example. \\
Our goal for the remaining of this work is to prove a relative version of this monoidal Bondal-Orlov theorem. We have to mention that such a theorem has been proven in the classical sense by Calabrese in \cite{calabrese2018relative}. Calabrese's precise formulation in his language is:
\begin{thm}\cite[Theorem 6.2]{calabrese2018relative}\label{thm:calabrese}
Let $S$ be a noetherian Artin stack with affine diagonal and let $X,Y\to S$ be flat, proper and relative algebraic spaces. Assume also that for al $s\in S$ the fibers $X_{s}, Y_{s}$ are projective, connected and Gorenstein and that $X_{s}$ has either ample canonical or anti-canonical bundle. Then $X\cong Y$ as $S$-stacks if and only if there exists a $S$-linear Fourier-Mukai equivalence $D(X)$.
\end{thm}
While there are many important aspects in the proof of this theorem, we are of the opinion that the key aspect in Calabrese's argument is the passing from fibers to neighborhoods in the Bondal-Orlov reconstruction argument. Concretely he shows
\begin{thm}\label{thm:calabresefiberstoglobal}
Let $S$ be a noetherian and connected scheme. Let $X$ be a noetherian Artin stack and $f:X\to S$ be a flat map. Let $E^{\bullet}$ be a bounded above complex with coherent cohomology. For any point $s\in S$, denote by $i_{\ast}:X_{s}\to X$ the inclusion of the fiber. Assume that for all $s\in S$ the derived restriction $i_{x}E^{\bullet}$ is concentrated in a single degree. Then $E^{\bullet}$ is a possibly shifted sheaf flat over $S$.
\end{thm}
It is this fibers-to-global result which will take the most content in this work to state and prove for our monoidal setting and which is for us too a key step in our main results. 
\section{From X to fibers}\label{sec:xtofibers}
This section will be devoted to the development of the necessary conditions our tensor triangulated categories will need to possess in order to be compatible with the fiber product of varieties. This is necessary as we start from an arbitrary tensor triangualted category structure $(\boxtimes,\mathbb{1})$ on $X\to S$ and we need to obtain an induced tensor triangualted structure on the derived category of fiber products $X_{Y}:=Y\times_{S}X$ for some other $S$-variety. \\
When working only with the underlying derived category or with this derived category equipped with the derived tensor structure induced by $Y$ this is of no issue since the fiber products $X_{Y}$ already have a variety structure and thus one can speak of the derived category without ambiguity even when base change theorems are not available for derived categories at the triangulated level. This is not true in general when one is concerned with the extra data of a TTC over the global derived category of $X\to S$.
We start with some basic definitions to settle some language:
\begin{defn}\label{def:transport}
Let $F\dashv G$ be an adjoint pair of exact functors $F:\T\to \T'$, $G:\T'\to \T$ between triangulated categories $\T,\T'$. Suppose $\T'$ is equipped with a TTC structure $(\boxtimes,\mathbb{1})$. We call the bifunctor
\[ F(G(\_)\boxtimes G(\_))\]
the transport of $\boxtimes$ along the adjoint pair $F\dashv G$. We denote this as $\boxtimes_{FG}$.
\end{defn}
We can analogously define a transport $_{GF}\boxtimes$ of a TTC structure $(\boxtimes,\mathbb{1})$ on $\T$ along the pair $F\dashv G$ although in our case we don't deal with structures transported in this direction. \\
There is an immediate problem with this definition. The transported structure $\boxtimes_{FG}$ does not generally traduce to a tensor triangulated category structure. \\
However, when the pair $F\dashv G$ comes from an equivalence of categories we do obtain a tensor structure when transporting along this equivalence. \\
Our concrete focus is on the following geometric situation: Let $\pi:X\to S=\mr{Spec} B$ be a projective smooth and faithfully flat morphism. Suppose that $(\boxtimes,\mathbb{1})$ is a TTC structure on $\der[x]$ which is $S$-linear in the sense that the action induced by the structure morphism $\pi$  is compatible with the tensor structure $(\dteee_{S},\Ox_{S})$ on $\der[S]$. In other words we want to have an action
\[ \der[S]\odot \der[X]\to \der[X] \]
Given by 
\[ K^{\bullet}\odot X^{\bullet}\mapsto \pi^{\ast}K^{\bullet}\boxtimes X^{\bullet} \]
Which additionally satisfies a projection formula
\[\pi_{\ast}(K^{\bullet}\odot X^{\bullet})\simeq K^{\bullet}\dteee_{S}\pi_{\ast}X^{\bullet}\]
Let us denote by $I_{s}:X_{s}\hookrightarrow X$ the inclusion of the fiber along a geometric $s:=\mr{Spec k}\to S$ (which we will in the future denote simply as $s\in S$). We would like to understand the case of tensor triangulated category structures $(\boxtimes,\mathbb{1})$ transported along the adjunction $I_{s}^{\ast}\dashv I_{s}\,_{\ast}$ between the bounded derived categories $\der[X_{s}]$ and $\der[X]$. \\
One immediate issue to deal with is that in general the transported bifunctor $\boxtimes_{FG}$ cannot be a monoidal structure simply because $F$ is not required to be a $\boxtimes$-functor and so we don't have 
\[ F(GX\boxtimes GX)\simeq FGX\boxtimes FGY.\]
And so the transported structure has no chance of being "associative" since 
\begin{align*}
(X\boxtimes_{FG}Y)\boxtimes_{FG} Z := \\F(G(X\boxtimes_{FG}Y)\boxtimes GZ)\simeq F(G(F(GX\boxtimes GY))\boxtimes GZ) \not\simeq \\ F(G((FGX\boxtimes FGY)\boxtimes FGZ)) \simeq X\boxtimes_{FG} (Y\boxtimes_{FG}Z)
\end{align*}
It can be shown however (\cite[Theorem 3.1]{kelly2006doctrinal}) that the transported structure does form an associative tensor triangulated category structure whenever the morphism
\[ F(X\boxtimes Y) \to F(GFX\boxtimes GFY)\]
induced by the counits, is an isomorphism of every pair of objects $X$ and $Y$. \\
As this is not generally true we are then forced to work with the non-associative transported structures. Nonetheless we can still phrase usual definitions from the theory of monoidal categories. 
\begin{defn}\label{def:unittransported}
    Let $F\dashv G$ be an adjoint pair of exact functors between triangulated categories $\T,\T'$. Suppose that $(\boxtimes,\mathbb{1})$ is a TTC structure on $\T$ and consider the transport $\boxtimes_{FG}$. We say that an object $U\in \T$ is a unit for $\boxtimes_{FG}$ if for any $X\in \T$ we have a natural equivalence
    \[ U\boxtimes_{FG}X \simeq  X\]
\end{defn}
The next definition is then natural:
\begin{defn}\label{def:invertible}
    Suppose $F\dashv G$, $\T,\T'$ and $(\boxtimes,\mathbb{1})$ be as above. Let $U\in \T$ be a $\boxtimes_{FG}$-unit. We say that $X\in \T$ is $\boxtimes_{FG}$-invertible if there exists $X^{-1}\in \T$ such that there are natural isomorphisms
    \[ X\boxtimes_{FG}X^{-1}\simeq U\]
\end{defn}
Let us remark here that $\boxtimes_{FG}$-invertible objects from Definition \ref{def:invertible} do not induce in general autoequivalences of the triangulated category $\T$. The fact that this happens when one deals with TTC structures is due to the associativity condition of monoidal categories, as otherwise these objects cannot induce equivalences of categories. \\
For us this means we have to make a sacrifice and ask for the better next thing from our structure $\boxtimes_{FG}$. We recall some basic definitions of the theory of monoidal categories first. 
\begin{defn}\label{def:rigid}
We say that a symmetric monoidal category $(\T,\boxtimes,\mathbb{1})$ is rigid if every object has a dual. Recall that for an object $X\in \T$, we say that a dual of $X$ is an object $X^{\vee}$ together with an evaluation map 
\[ X^{\vee}\otimes X \to \mathbb{1} \]
Satisfying natural compatibility conditions with the monoidal structure.
\end{defn}
As a consequence, a rigid symmetric monoidal category is automatically closed and for every object $X$ the functors $X\boxtimes\_$ have left and right adjoint functors given by $X^{\vee}\boxtimes\_$. \\
This condition is relatively mild in applications in algebraic geometry, it is indeed true that derived categories of perfect complexes with monoidal structure given by $\dteee$ are rigid under the usual definition of dual object. \\
\\
An important invariant of derived categories is the Serre functor as defined by Kapranov and Bondal, let us recall:
\begin{defn}\label{def:serrefunctor}
Let $\T$ be a triangulated category, a functor $\mathbb{S}:\T\to\T$ is a Serre functor if we have isomorphisms:
\[ \mr{Hom}_{\T}(A,B)\cong \mr{Hom}_{\T}(B,\mathbb{S}(A))^{\ast}\]
For every $A,B\in \T$
\end{defn}
Kapranov and Bondal showed that whenever a Serre functor exists, it is unique up to a graded equivalence. In $\der[X]$ for $X$ a smooth projective variety of dimension $n$, the Serre functor can be shown to be equivalent to $\_\dtee\omega_{X}[n]$.  \\
One important property of Serre functors is that one can see that they commmute with autoequivalences. In other words, if $\phi:\T\to \T$ is in $\mr{Aut}(\T)$ then $\phi\circ \mathbb{S}\simeq \mathbb{S}\circ \phi$. \\
this result was a very important ingredient in the proof of the monoidal Bondal-Orlov theorem in \cite{toledo2024tensor}, as we kept some control in the tensor structure $\boxtimes$ through $\boxtimes$-invertible objects since they induce autoequivalences and therefore must interact nicely with the functor $\dteee\omega_{X}[n]$. \\
As discussed above, since $\boxtimes_{FG}$-invertible objects don't necessarily induce autoequivalences we cannot use this property of Serre functors in our current relative situation. \\
There is still some compatbility with Serre functors even in the case when we are dealing with non-autoequivalences, for example:
\begin{thm}\label{thm:serreadj}
    Let $F:\A\to \B$ be a functor between $k$-linear categories. Suppose both $\A$ and $\B$ are endowed with Serre functors $\mathbb{S}_{\A}$ and $\mathbb{S}_{\B}$ respectively. Then 
    \[ G\dashv F \Rightarrow F\dashv \mathbb{S}_{\A}\circ G\circ \mathbb{S}_{\B}^{-1}.\]
    \begin{proof}
        \cite[Remark 1.31]{huybrechts2006fourier}.
    \end{proof}
\end{thm}
In other words, we can calculate right adjoints as $\mathbb{S}_{\A}\circ G\mathbb{S}^{-1}_{\B}$. \\
Our goal here is to be able to compare the TTC structure given by the derived tensor product $\dteee_{X_{s}}$ with the structure of a transported TTC structure on $\der[X]$ to $\der[X_{s}]$. \\
\begin{defn}\label{def:relduals}
    Suppose $(\boxtimes,\mathbb{1})$ is a TTC structure on a triangulated category $\T$ and we have an adjunction $F\dashv G$. We say that an object $X$ has a dual $X^{\vee}$ with respect to $\boxtimes_{FG}$ if there is a morphism, antural in $X^{\vee}$
    \[ X^{\vee}\boxtimes_{FG} X\to \mathbb{1}\]
    Such that these induce isomorphisms
    \[ \mr{Hom}_{\T}(A,B\boxtimes_{FG}X^{\vee}) \to \mr{Hom}_{\T}(A\boxtimes_{FG} X, B)\]
    For any pair of objects $A,B\in \T$. 
\end{defn}
We will furthermore work under the assumption that all these dual objects are strong duals in the sense that the natural morphisms
\[ X\mapsto X^{\vee\vee} \]
are isomorphisms. \\
One can show under these hypothesis, and using that the product $\boxtimes_{FG}$ is symmetric, that we can obtain adjunctions
\[ X^{\vee}\boxtimes_{FG} \_ \dashv X\boxtimes_{FG}\_\dashv X^{\vee}\boxtimes_{FG}\_ \]
To summarize our assumptions here, we introduce the following definition:
\begin{defn}\label{def:reasonable}
    Let $F\dashv G:\T\to \T'$ be an adjoint pair between triangulated categories $\T, \T'$. We say that a tensor triangulated category structure $(\boxtimes,\mathbb{1})$ on $\T'$ is reasonable with respect to $F\dashv G$ if the transported structure $\boxtimes_{FG}$ satisfies the following conditions:
    \begin{enumerate}
        \item The object $F(\mathbb{1})$ is a unit object for $\boxtimes_{FG}$. 
        \item Every invertible object has strong duals.
        \item If $X$ is $\boxtimes_{FG}$-invertible, then the functor $X\boxtimes_{FG}\_$ is full. 
    \end{enumerate}
\end{defn}
Let us now fix a TTC structure $(\boxtimes,\mathbb{1})$ in $\der[X]$ which we will assume from now on, is reasonable with respect to the pair $I_{s}^{\ast}\dashv I_{s}\,_{\ast}:\der[X]\to\der[X_{s}]$ denoting the derived pullback and pushforward functors induced by the inclusion of the fiber $I_{s}:X_{s}\hookrightarrow X$ along a point $s\in S$. 
\begin{lemma}\label{lemma:serretransported}
Let $(\boxtimes, \mathbb{1})$ be a reasonable TTC structure on $\der[X]$ and let $A\in \der[X_{s}]$ be a $\boxtimes_{I_{s}^{\ast}I_{s}\,_{\ast}}$-invertible object. Then we have
\begin{align}\label{eq:serretransported} A^{\vee}\boxtimes_{I_{s}^{\ast}I_{s}\,_{\ast}}(\omega_{X_{s}}\dteee_{X_{s}} \_) \simeq \omega_{X_{s}}\dteee_{X_{s}}(A^{\vee}\boxtimes_{I_{s}^{\ast}I_{s}\,_{\ast}}\_) \end{align}
\begin{proof}
    Immediate from Theorem \ref{thm:serreadj}, that the Serre functor in $\der[X_{s}]$ is given by $\omega_{X_{s}}\dteee_{X_{s}}\_[n_{s}]$ (where $n_{s}$ denotes the dimension of $X_{s}$), and the adjunction 
    \[ A^{\vee}\boxtimes_{I_{s}^{\ast}I_{s}\,_{\ast}} \dashv A\boxtimes_{I_{s}^{\ast}I_{s}\,_{\ast}} \dashv A^{\vee}\boxtimes_{I_{s}^{\ast}I_{s}\,_{\ast}}.\]
\end{proof}
\end{lemma}
\begin{lemma}\label{lemma:invdualinv}
    Suppose $(\boxtimes,\mathbb{1})$ is a reasonable TTC structure with respect to an adjoint pair $F\dashv G:\der[X]\to \der[Y]$. Suppose $A$ is a $\boxtimes_{FG}$-invertible object, then $A^{\vee}$ is also $\boxtimes_{FG}$-invertible. 
    \begin{proof}
        The proof follows a standard argument. By Yoneda, the isomorphism
        \[y:=A\boxtimes_{FG}A^{-1}\to V:=F(\mathbb{1})\]
        induces a natural transformation 
        \[ \mr{Hom}_{\der[X]}(\_,A^{-1})\to \mr{Hom}_{\der[X]}(\_\boxtimes_{FG}A,V).\]
        Which by the fact that $\_\boxtimes_{FG}A$ is full, and thus faithful, induces a bijection
        \[ \mr{Hom}_{\der[X]}(A^{\vee},A^{-1})\to \mr{Hom}_{\der[X]}(A^{\vee}\boxtimes_{FG}A, V).\]
        If we consider the evaluation $\epsilon:A^{\vee}\boxtimes_{FG}A\to V$, it corresponds by the previous equivalence to the morphism $\phi:A^{\vee}\to A^{-1}$ such that \begin{align*}
        \xymatrix{\epsilon=A^{\vee}\boxtimes_{FG}A \ar[rr]^{\phi\boxtimes_{FG}1_{A}} && A^{-1}\boxtimes_{FG} A\ar[r]^{y} & V.}\end{align*}
        This completes the proof
        \end{proof}
\end{lemma}
We now move on to more our more specific geometric situation. 
\begin{lemma}\label{lemma:compttcfibers}
Let $s\in S$ and $\boxtimes_{I_{s}^{\ast}I_{s}\,_{\ast}}$ be as above and suppose that the canonical bundle $\omega_{X_{s}}$ of $X_{s}$ is $\boxtimes_{I_{s}^{\ast}I_{s}\,_{\ast}}$-invertible. Then we have
\begin{align}\label{eq:compttcfibers}
    \omega_{X_{s}}\boxtimes_{I_{s}^{\ast}I_{s}\,_{\ast}}(\omega_{X_{s}}\dteee_{X} \_ )\simeq \omega_{X_{s}}\dtee(\omega_{X_{s}}\boxtimes_{I_{s}^{\ast}I_{s}\,_{\ast}}\_).
\end{align}
\begin{proof}
    As we have an isomorphism $\omega^{\vee\vee}_{X_{s}}\cong \omega_{X_{s}}$ by Lemma \ref{lemma:invdualinv}, if we consider the $\boxtimes_{I_{s}^{\ast}I_{s}\,_{\ast}}$-invertible object $\omega_{X_{s}}$ we obtain our result by use of Lemma \ref{lemma:serretransported}.
\end{proof}
\end{lemma}
We immediately get :
\begin{cor}\label{cor:cbboxcb}
    Let $\boxtimes_{I_{s}^{\ast}I_{s}\,_{\ast}}$ be as above and suppose $I^{\ast}(\mathbb{1})\simeq \Ox_{X_{s}}$. Furthermore, suppose $\Ox_{X_{s}}$ is a unit for $\boxtimes_{I_{s}^{\ast}I_{s}\,_{\ast}}$. Then we have
    \begin{align}\label{eq:cbboxcb}
        \omega_{X_{s}}\boxtimes_{I_{s}^{\ast}I_{s}\,_{\ast}} \omega_{X_{s}} \cong \omega_{X_{s}}\dteee_{X_{s}}\omega_{X_{s}}.
    \end{align}
\end{cor}
It is this latter identity which will allow us to fully compare these transported structures $\boxtimes_{I_{s}^{\ast}I_{s}\,_{\ast}}$ with the derived tensor product $\dteee_{X_{s}}$ of $\der[X_{s}]$.
\begin{thm}\label{thm:mainthm1}
Let $\pi:X\to S=\mr{Spec} B$ be a smooth projective variety faithfully flat over $S$. Let $(\boxtimes,\mathbb{1})$ be a $S$-linear tensor triangulated category structure on $\der[X]$. Let $s\in S$ be a closed point and let $I_{s}:X_{s}\hookrightarrow X$ denote the inclusion of the fiber $X_{s}$ which we suppose has either ample canonical or anti-canonical bundle $\omega_{X_{s}}$. Suppose $I_{s}^{\ast}\mathbb{1}\cong \Ox_{X_{s}}\in \der[X_{s}]$ and that $(\boxtimes,\mathbb{1})$ is reasonable with respect to the adjunction $I_{s}^{\ast}\dashv I_{s}\,_{\ast}:\der[X]\to \der[X_{s}]$. If $\omega_{X_{s}}$ is $\boxtimes_{I_{s}^{\ast}I_{s}\,_{\ast}}$-invertible then $\boxtimes_{I_{s}^{\ast}I_{s}\,_{\ast}}$ coincides on objects with $\dteee_{X_{s}}$.
\begin{proof}
    Since the fiber $X_{s}$ has an ample canonical bundle $\omega_{X_{s}}$, the collection $\{\omega_{X_{s}}^{\dteee i}\}_{i}$ forms an ample sequence and thus a spanning class of $\der[X_{s}]$. this in turn means that 
    \[ \bigoplus_{i,j}\omega_{X_{s}}^{\dteee i}[j] \]
    is a weak generator and therefore a strong generator of $\der[X_{s}]$. \\
    In \cite[Theorem 3.10]{toledo2024tensor} the second author shows that this fact can be used to conclude that any object $A^{\bullet}\in \der[X_{s}]$ is generated by shifts of direct sums of objects of the form $\omega_{X_{s}}^{\dteee i}$ for some $i\in \mathbb{Z}$. \\
    Since $\boxtimes_{I_{s}^{\ast}I_{s}\,_{\ast}}$ is exact on both variables, we use Lemma \ref{cor:cbboxcb} to conclude that 
    \[ \omega_{X_{s}}\boxtimes_{I_{s}^{\ast}I_{s}\,_{\ast}} A^{\bullet} \simeq \omega_{X_{s}}\dteee_{X_{s}} A^{\bullet}   \]
    for any $A^{\bullet}\in \der[X]$. \\
    Similarly, for any $B^{\bullet}$ we have
    \[ B^{\bullet}\boxtimes_{I_{s}^{\ast}I_{s}\,_{\ast}} A^{\bullet} \simeq B^{\bullet}\dteee_{X_{s}} A^{\bullet}.\]
\end{proof}
\end{thm}
Before continuing, let us briefly remark that under the hypothesis of the proposition above one could conclude that the transported structure $\boxtimes_{I_{s}^{\ast}I_{s}\,_{\ast}}$ is associative as it agrees pairwise with the monoidal category structure given by $\dteee_{X_{s}}$. While this is true in the sense that we obtain isomorphisms 
\[ (A^{\bullet}\boxtimes_{I_{s}^{\ast}I_{s}\,_{\ast}} B^{\bullet})\boxtimes_{I_{s}^{\ast}I_{s}\,_{\ast}} C^{\bullet}\cong A^{\bullet}\boxtimes_{I_{s}^{\ast}I_{s}\,_{\ast}} (B^{\bullet}\boxtimes_{I_{s}^{\ast}I_{s}\,_{\ast}} C^{\bullet}). \]
This conclusions is however \textit{a fortiori} and doesn't say anything of the associative coherent conditions of a monoidal category or even the associator morphisms. \\
This is of no issue when one deals with Balmer spectra since the construction of $\mr{Spec}(\boxtimes)$ actually does not depend on the particular choice of associators nor on whether they satisfy any coherence condition since the construction only uses the data of the bifunctor. We will come back to this point in the next section. 
\begin{obs}\label{obs:magmoidal}
In the literature, the data of a bifunctor $\boxtimes:C\times C\to C$ for some category $C$ and with a unit condition is sometimes refered as a magmoidal category. In this sense one could extend Balmer's construction of to what we ought to call magmoidal triangulated categories. As to our knowledge there are no immediate reasons to consider that is either substantially different or which can be applied to a situation not covered by the usual considerations when studying Balmer spectra, we do not develop this idea in further generality although our transported structures $\boxtimes_{FG}$ could already be called magmoidal structures.
\end{obs}
\section{Enhancing tensor structures}\label{sec:ttenhancement}
One classical subject of interest in the study of derived categories of sheaves in algebraic geometry (and generally in the study of triangulated categories)  is that of the existence of higher enhancements. As a great deal of literature has been dedicated to the subject we will only cover the strictly necessary results in this section which will be of use later during this work. \\
In this section we first give a very brief overview of this theory and of Toën's Morita theorem for dg-categories including the description of the Dwyer-Kan model structure on the category of dg-categories and the concrete applications we have in mind in algebraic geometry. \\
Next, we recall some of the constructions developed in \cite{toledo2023davydov} which describe enhancements of magmoidal triangulated category structures. 
\subsection{The theory of dg-enhancements}
While enhancements of triangulated categories can be performed using different languages such as derivators, stable $\infty$-categories and general $A_{\infty}$-categories, we pick dg-categories as our working setting. Let us recall some basics of the theory of such objects. For an interested reader looking for a more in depth treatment of the use of dg-categories in algebraic geometry we point towards \cite{toen2008lectures} and \cite{canonaco2017tour}.
\begin{defn}\label{def:dgcategory}
    A dg-category over a ring $R$ is a category enriched over the category of chain complexes $\cc$. 
\end{defn}
Following the classical theory of enriched category theory, we have enriched versions of the many categorical constructions. In particular of enriched functors and natural transformations. We shall call these, respectively, dg-functors and dg-natural transformations. 
\begin{defn}\label{def:dgmodules}
    Let $T$ be a dg-category over $R$, we will call the category of dg-functors $\mr{Fun}(T^{op},\cc)$ as the category of right dg-modules over $T$ and denote it by $\mr{T^{op}-Mod}$. 
\end{defn}
There is also an enriched version of the classical Yoneda embedding.
\begin{defn}\label{def:enrichedyoneda}
    Let $T$ be a dg-category. The Yoneda embedding is the dg-functor
    \[ h_{\_}:T\to T^{op}-\mr{Mod}\]
    Which maps $x\in T$ to $h_{x}:=Hom_{T}(\_,x)$.
\end{defn}
As usual, we say that a left $T$-module is representable if it is equivalent as a dg-functor in $T^{op}$-Mod to a right module of the form $h_{x}$ for some $x\in T$. \\
To a given dg-category $T$ one can associate a category $H^{0}(T)$ defined by having the same objects as $T$ and for any pair of objects $x,y\in H^{0}(T)$ we define 
\[Hom_{H^{0}(T)}(x,y):=H^{0}(Hom(x,y)) \]
Similarly a dg-functor $F:T\to T'$ induces a functor $H^{0}(F):H^{0}(T)\to H^{0}(T')$. \\
As we have mentioned before, our use for dg-categories all through this work is centered around the idea of enhancing triangulated categories. Roughly speaking, the idea behind this is that for a given triangulated category $\T$, one finds a dg-category $T$ which holds all information of $\T$ in the homotopy category $H^{0}(T)$ plus important homotopical information of it. \\
In this sense, the category $\T$ is a ``shadow" of the dg-category $\T$. We now give a series of formal definitions related to this idea. 
\begin{defn}\label{def:dgenhancement}
    Let $\T$ be a triangulated category. We say that $T$ is a dg-enhancement of $\T$ if there exists a triangulated category equivalence
    \[ \epsilon:H^{0}(T)\to \T. \]
\end{defn}
It can be shown that the category $H^{0}(\mr{T^{op}-Mod})$ can be equipped with a natural triangulated category structure. A consequence of this serves as a justification for the following definition.
\begin{defn}\label{def:pretriangulated}
    Let $T$ be a dg-category. We say that it is pretriangulated if the image of the Yoneda embedding functor 
    \[H^{0}(h_{\_}):H^{0}(T)\to H^{0}(\mr{T^{op}-Mod}) \]
    is a triangulated subcategory.
\end{defn}
\begin{defn}\label{def:triangulatedhull}
    If $T$ is a dg-category, we write $T^{\mr{pre-tr}}$ for the smallest pretriangulated full dg-subcategory of $\mr{T^{op}-Mod}$. 
\end{defn}
We can think of the process of applying the $(\_)^{\mr{pre-tr}}$ operation as formally adding mapping cones, direct sums, etc. \\
Let us also write $\mr{tri}(T)$ for the triangulated category $H^{0}(T^{\mr{pre-tr}})$ and $Perf(T)$ for the full subcategory of compact objects in $\mr{tri}(T)$.\\
\begin{defn}\label{def:dgquotient}
    Let $T$ be a dg category and let $S\subset T$ be a full sub dg-category. The quotient $T/S$ is defined as the dg-category with the same collection of objects as $T$ and such that for every $s\in S$, we add an endomorphism $\eta$ of degree $-1$ so that $d(\eta)=Id_{s}$.
\end{defn}
With this definition of the quotient of dg-categories, it is possible to show that $H^{0}(T)/H^{0}(S)\simeq H^{0}(T/S)$ and so the derived category of an abelian category $\A$ can be enhanced by the quotient $\cc[\A]/Ac(\A)$, where $Ac(\A)$ is the subcategory of acyclic complexes of $\A$. \\
We can similarly argue about the bounded (below, above, or both) derived categories by simply taking the corresponding subcategories of $\cc[\A]$ and $Ac(\A)$. \\
\begin{defn}\label{def:dgderivedcat}
Let $T$ be a dg-category, write $Ac(T)$ for the full dg-subcategory of $T^{op}-Mod$ spanned by those modules which are object-wise acyclic. The derived category $D(T)$ of $T$ is the quotient $H^{0}(T^{op}-\mr{Mod}/Ac(T))$.
\end{defn}
This derived category is always triangulated as it is a Verdier quotient of the triangulated category $H^{0}(T^{op}-\mr{Mod})$. \\
This is in particular implies that for a scheme there always exists a dg-enhancement of its derived category. \\
The next natural question is whether a functor between triangulated categories can be lifted to a functor between dg-enhancements of the domain and target. As it turns out, in order for this idea to work we need to study the formal homotopy theory of dg-categories. \\
In \cite{tabuada2007theorie}, Tabuada described a model category structure on the category of dg-categories given as a Dwyer-Kan model category in the context of dg-categories. We describe the weak equivalences and fibrations.
\begin{defn}\label{def:DKdgmodel}
    Let $T, T'$ be dg-categories over $R$. We say that a dg-functor $F:T\to T'$ is a weak equivalence if
    \begin{enumerate}
        \item It is quasi-fully faithful. This means that the induced cochain complex morphism $Hom_{T}(x,y)\to Hom_{T'}(F(x),F(y))$ is a quasi-isomorphism for all $x,y\in T$. 
        \item It is quasi-essentially surjective: For any $x'\in H^{0}(T')$, there exists $x\in T$ such that $H^{0}(F)(x)\simeq x'$. 
    \end{enumerate}
    Fibrations are in turn described as follows
    \begin{enumerate}
        \item The induced morphism of complexes $Hom_{T}(x,y)\to Hom_{T'}(F(x),F(y))$ is a fibration in the (projective) model structure of unbounded complexes.
        \item For any isomorphism $u':x'\to y'\in H^{0}(T)$ and any $y\in H^{0}(T)$ such that $F(y)=y'$ there is an isomorphism $u:x\to y$ in $H^{0}(T)$ such that $H^{0}(F)(u)=u'$.
    \end{enumerate}
\end{defn}
These classes of morphisms in the category of dg-categories forms a model category. The resulting homotopy category of this structure will be denoted by $H_{qe}$. \\
We say that a dg-functor $f\in H_{qe}(T,T')$ is a quasi-functor. \\
With this, we can describe enhancements of triangulated functors. \\
\begin{defn}\label{def:dglifts}
    Let $\T,\T'$ be triangulated categories with dg-enhancements $T, T'$. We say that an exact functor $\F:\T\to \F$ has a dg-lift if there exists a morphism $f\in H_{qe}(T,T')$ such that $H^{0}(f)=\F$. 
\end{defn}
We can now give a notion of uniqueness of dg-enhancements.
\begin{defn}\label{def:uniquedg}
    Let $\T$ be a triangulated category and let $T$, $S$ be two enhancements $\epsilon:H^{0}(T)\to \T$, $\epsilon':H^{0}(S)\to \T$. We say that $\T$ has a unique enhancement if there exists a quasi-functor $f:T\to S$ such that $H^{0}(f)$ is an equivalence of triangulated categories. 
\end{defn}
From now on, let us fix and denote by $Perf^{dg}(X)$ a dg-enhancement of the derived category of perfect complexes over a space (variety, scheme, etc) $X$, and similarly we put $QCoh^{dg}(X)$ for a dg-enhancement of the derived category of quasi-coherent sheaves on $X$. \\
A fundamental result in the theory of dg-categories and their homotopy theory is Toën's Morita theorem proven in \cite{toen2007homotopy}. In it, Toën shows that while the enriched tensor product of dg-categories is not compatible with the model category structure described above, it is nonetheless possible to equip $H_{qe}$ with a closed symmetric monoidal structure by considering a derived tensor product $\dteee_{dg}$ from which it is possible to describe the internal Hom functor $\ihom{\_}{\_}$. \\
Given a model category $M$ we will write $Int(M)$ for the full subcategory of fibrant and cofibrant objects of $M$. In our case we are interested in dg-categories of modules $T^{op}-\mr{Mod}$ over a dg-category $T$, this category can be equipped with a model category structure: \\
\begin{enumerate}\label{def:modeldgmod}
    \item Weak morphisms $f:F\to G\in T^{op}-\mr{Mod}$ are those such that for any $x\in T^{op}$, the induced morphism $f_{x}:F(x)\to G(x)$ is a quasi-isomorphism.
    \item Fibrations are those $f:F\to G$ such that the morphism $f_{x}$ is a fibration in $\cc$.
\end{enumerate}
 Importantly, the dg-modules $h_{x}:=Hom_{T}(\_,x)$ are both fibrant and cofibrant and so the Yoneda embedding of Definition \ref{def:enrichedyoneda} factorizes through $Int(T^{op}-\mr{Mod})$. 
 \begin{defn}\label{def:rqr}
     Let $T$ be a dg-category. A dg-module $F\in T^{op}-\mr{Mod}$ is called quasi-representable if it is equivalent in $Ho(T^{op}-\mr{Mod})$ to a module of the form $h_{x}$. 
 \end{defn}
We denote by $T^{op}-\mr{Mod}^{rqr}$ the category of all quasi-representable right modules over $T$. \\
\begin{defn}\label{def:bimodrqr}
For a collection of dg-categories $\{T_{1},\dots,T_{k}\}$, and a dg-category $S$, we say that a dg-module 
\[F \in T_{1}\dteee_{dg}\dots\dteee_{dg} T_{k}\dteee_{dg}S^{op}-\mr{Mod}\] 
is right quasi-representable if for any collection of objects $x_{i}\in T_{i}$ the dg-functor \[F(x_{1},\dots,x_{k}):S^{op}\to \cc\] 
is right quasi-representable. \\
We denote by $T_{1}\dteee_{dg}\dots\dteee T_{k}\dteee_{dg} S^{op}-\mr{Mod}^{rqr}$ the dg-category of right quasi-representable $T_{1}\dteee_{dg}\dots\dteee_{dg} T_{k}\dteee_{dg}S^{op}$-modules.
\end{defn}
Toën provides a characterization of the internal derived $Hom$ $\ihom{T}{T'}$ for $T,T'$ in $H_{qe}$ as
\begin{align}\label{def:internalhom}
\ihom{T}{T'}\simeq Int(T\dteee_{dg}T'^{op}-\mr{Mod}^{rqr}).
\end{align}
In particular we have that
\[ \hat{T}:=Int(T^{op}-\mr{Mod})\simeq \ihom{T^{op}}{Int(\cc)}.\]
We write $T_{pe}$ for the dg-category of compact objects of $\hat{T}$. \\
Let us denote by $\ihomc{T}{T'}$ the category of continuous dg-functors between dg-categories $T, T'$. \\
We can now formulate the main theorem of \cite{toen2007homotopy}.
\begin{thm}\label{thm:dgmorita}
    Let $T,S$ be dg-categories and let $h:T\to T^{op}-\mr{Mod}$ denote the Yoneda embedding. For any other dg-category $S$, we have
    \begin{enumerate}
        \item The pullback functor $h^{\ast}:\ihomc{\hat{T}}{\hat{S}}\to \ihom{T}{\hat{S}}$ is an isomorphism in $H_{qe}$.
        \item The pullback functor $h^{\ast}:\ihom{T_{pe}}{S_{pe}}\to \ihom{T}{S_{pe}}$ is an equivalence in $H_{qe}$.
    \end{enumerate}
\end{thm}
From this a more familiar form of the result follows:
\begin{cor}\label{cor:dgmorita}
    Let $T$ and $S$ be two dg-categories, then there exists a natural isomorphism $H_{qe}$
    \[ \ihomc{\hat{T}}{\hat{S}}\simeq \widehat{T^{op}\dteee_{dg}S}.\]
\end{cor}
A consequence of Theorem \ref{thm:dgmorita} in the geometric context is 
\begin{thm}\label{thm:eilenbergwatts}
    Let $X,Y$ be two quasi-compact quasi-separated schemes over $R$, assume one is flat over $R$. Then there exists an isomorphism in $H_{qe}$:
    \[ \ihomc{QCoh^{dg}(X)}{QCoh^{dg}(Y)}\simeq QCoh^{dg}(X\times_{R} Y).\]
\end{thm}
And finally, 
\begin{thm}\label{thm:perfecteilenbergwatts}
    Let $X,Y$ be two smooth and proper schemes over a field $k$. Then there exists an isomorphism in $H_{qe}$:
    \[ \ihom{Perf^{dg}(X)}{Perf^{dg}(Y)}\simeq Perf^{dg}(X\times_{k} Y).\]
\end{thm}
This result establishes a correspondence between dg-lifts of exact functors between derived categories of perfect complexes and objects in $Perf^{dg}(X\times_{k} Y)$ though of as kernels of Fourier-Mukai transforms.
\subsection{Lifting tensor products}
For us, the importance of Theorem \ref{thm:perfecteilenbergwatts} relies in the fact that we can encode the information of bifunctors $\boxtimes:\T\times\T\to \T$ in terms of certain dg-bimodules at the level of dg-enhancements of $T$. \\
In \cite{toledo2023davydov}, the second author carried an exploration of this situation and settled on a truncated version of dg-enhancements of tensor triangulated category structures. In our particular situation, as pointed out in Remark \ref{obs:magmoidal} we are interested in enhancing the data of the bifunctor and we don't necessarily have any coherence data to worry about. \\
We present a version of the work as in \cite{toledo2023davydov} limited only to the case of Magmoidal Triangulated Categories (MTC's for short). 
\begin{defn}\label{def:dgbimod}
    Let $T$ be a dg-category. An n-fold dg-bimodule over $T$ is a dg-module $F\in T^{\tofdgc n}\tofdgc T^{op}-\mr{Mod}$.
\end{defn}
This is a slight generalization of dg-bimodules over a dg-category $T$, which we can identify with 1-fold dg-bimodules. There is a dg-category $\mr{Bimod}^{n}_{dg}(T)$ of n-fold dg-bimodules with morphisms given by dg-natural transformations. \\
One can also take tensor products of an n-fold dg-bimodule $\Gamma$, and an m-fold dg-bimodule $\Lambda$. Concretely, if 
\[\Gamma:\underbrace{T\tofdgc\dots\tofdgc T}_{n} \tofdgc T^{op}\to \cc, \Lambda:\underbrace{T\tofdgc\dots\tofdgc T}_{m}\tofdgc T^{op}\to \cc \]
then we get
\[ \mr{Bimod}_{dg}^{n+m-1}(T)\ni \Gamma\otimes \Lambda : T^{\tofdgc n}\tofdgc T^{op}\tofdgc T \tofdgc T^{\tofdgc m-1}\tofdgc T^{op}\to \cc.\]
By which we mean we use the rightmost factor $T^{op}$ of $\Gamma$ and the leftmost factor $T$ of $\Lambda$ to construct this tensor product as one does with any tensor product of bimodules (see for example  \cite{drinfeld2004dg} for an explicit construction). \\
We borrow the following notation from \cite{hovey2011additive}: If we wanted to take this tensor product using any other of the left $T$ factors of $\Lambda$, we could make this explicitiy by putting labels on each of these factors so that for example if $\Gamma\in \mr{Bimod}^{2}_{dg}(T)$, we can write $\Gamma_{X,Y}$ where $X,Y$ are to be understood of labels for each of the left factors $T$ in $\Gamma:T\tofdgc T\tofdgc T^{op}\to \cc$. \\
Then, if we pick for example another 2-fold dg-bimodule $\Lambda$ and we want to calculate the 3-fold dg-bimodule $\Gamma\otimes \Lambda$ defined by the tensor product of the rightmost $T^{op}$ factor of $\Gamma$ and the first $T$ factor in $\Lambda$, we write this as 
\[ \Gamma_{X,Y}\otimes \Lambda_{Z,\Gamma}:\overset{\textbf{X}}{T}\tofdgc \overset{\textbf{Y}}{T}\tofdgc T^{op} \tofdgc \overset{Z}{\textbf{T}}\tofdgc \overset{\Gamma}{T} \tofdgc T^{op}\to \cc. \]
By which we mean that it is the dg-category $T$ labeled with $\Gamma$ on the right that will be tensored with the first appearing $T^{op}$ belonging to the dg-bimodule $\Gamma$. 
In \cite[Theorem 2.31]{toledo2023davydov} the second author proved the following theorem:
\begin{thm}\label{thm:eudglift}
    Let $A$ be a dg-algebra and let $\boxtimes:\hzero[A_{pe}]\times\hzero[A_{pe}]\to \hzero[A_{pe}]$ be an exact functor in each variable. Suppose that for every object $M\in \hzero[A_{pe}]$, the triangulated functors
\begin{align*}
    M\boxtimes \_ :\hzero[A_{pe}] \to \hzero[A_{pe}]
    \\
    \_\boxtimes M:\hzero[A_{pe}] \to \hzero[A_{pe}]
\end{align*}
both have unique dg-enhancements $R_{M}$ and $L_{M}$ respectively. \\
Then, $L_{A}(A)$ is a 2-fold dg-bimodule and for any $N\in A_{pe}$ we have
\[ \hzero[L_{A}(A)\otimes M\otimes N]\simeq M\boxtimes N. \]
\end{thm}
The hypothesis here is relevant to our situation as it is well known that the dg-enhancements $Pef^{dg}(X)$ of a smooth projective variety are all equivalent in $H_{qe}$ to the category of perfect complexes over some dg-algebra $A$. The condition on the existence of dg-lifts is not a strong one either as one can show that this condition is equivalent to the functors being of Fourier-Mukai type. The uniqueness condition is more delicate but in practice it is often satisfied.  \\
With this result in place, we define
\begin{defn}\label{def:pseudodgmagmoidal}
A (left-biased lax) pseudo dg-magmoidal structure in a dg-category $T$ consists on the data:
\begin{enumerate}
    \item A 2-fold right quasi-representable dg-bimodule $\Gamma \in \mr{Bimod}^{2}_{dg}(T)$.
    \item A quasi-representable n object $U\in T^{op}-\mr{Mod}$ called the unit.
    \item  Morphisms of 3-fold dg-bimodules $\alpha_{X,Y,Z}:\Gamma_{Y,Z}\otimes\Gamma_{X,\Gamma}\to  \Gamma_{X,Y}\otimes\Gamma_{\Gamma,Z}$.
    \item A  morphism of 1-fold dg-bimodules $\ell_{X}:U\otimes\Gamma_{U,X} \to h_{X}$.
    \item A  morphism of 1-fold dg-bimodules $r_{X}:U\otimes\Gamma_{X,U}\to h_{X}$.
    \item A morphism $c_{X,Y}:\Gamma_{X,Y} \to \Gamma_{Y,X}$ of 2-fold dg-bimodules.
\end{enumerate}
We require that the morphisms $\alpha_{X,Y,Z}$, $u_{X}$ and $c_{X,Y}$ are all isomorphisms when passing to the homotopy category $\hzero[\cc]$.  \\
Furthermore we require the existence of the following homotopy data satisfying conditions:
\begin{enumerate}
    \item (Unit) A morphism $\mu\in Hom^{-1}(U\otimes\Gamma_{X,\Gamma}\otimes \Gamma_{U,Y}, U\otimes\Gamma_{\Gamma,Y}\otimes \Gamma_{X,U}) $ such that  $\alpha_{X,U,Y}^{0}\circ Id_{X}^{0}\otimes \ell^{0}_{Y}-r_{X}^{0}\otimes Id_{Y}^{0}=d(\mu)$.
    \item (Symmetry) The composition $c_{X,Y}\circ c_{Y,X}$ is the identity in $\hzero[T]$. 
    \item (Unit symmetry) There is $\kappa\in Hom^{-1}(U\otimes \Gamma_{X,U}, h_{X})$ such that $\ell_{X}\circ c_{X,U}-r_{X}=d(\kappa)$
\end{enumerate}
\end{defn}
This definition is too general and for our applications we need to restrict to those pseudo dg-magmoidal structures which preserve perfect objects.
\begin{defn}\label{def:perfectdgtensor}
A pseudo dg-magmoidal structure in a dg-category $T$ with 2-fold dg-bimodule $\Gamma$ is called perfect if for every $X,Y\in T$, $X\otimes Y\otimes \Gamma_{X,Y}$ is quasi-represented by a perfect $T^{op}$-module.
\end{defn}
Finally we present the following theorem which can be proven just as in \cite[Lemma 2.26]{toledo2023davydov}
\begin{lemma}\label{lemma:correspondence}
A perfect pseudo dg-magmoidal structure $\Gamma$ on a dg-category $T$ induces a magmoidal triangulated category structure (see \ref{obs:magmoidal}) on $\hzero[T_{pe}]$.
\end{lemma}
As a consequence of this, if we consider a triangulated category $\T$ which has a dg-enhancement $T$, a perfect pseudo dg-magmoidal (pp dg-magmoidal for short) induces a MTC on $\T$.  \\
\section{Descent of dg-categories}\label{sec:descent}
In this section we quickly unwrap one of the main theorems of Hirschowitz and Simpson in \cite{hirschowitz1998descente} which is behind our indiscriminate use of two important tools in this work: That the assignment to each open $U\subseteq S$ of its derived category does glue into a geometric object, namely a stack, when one deals with dg-categories, and that it is enough to work over an affine base $S=\mr{Spec} B$. \\
It is these tools combined with some formal properties of the closed symmetric monoidal structure on the $\infty$-category of dg-categories as shown by Toën in \cite{toen2007homotopy} which will ultimately allow us to work with a stack of tensor triangulated categories. \\
We will start by giving a very quick summary of the notation and definitions important in the main theorem of Hirschowitz and Simpson, we refer to \cite{hirschowitz1998descente} for the full details of their framework.  \\
The following is a paraphrased translation from the french taken from their manuscript. \\
\begin{defn}\label{def:siteeds}\cite[Chapter 15]{hirschowitz1998descente}
We say that a site $\X$ admits enough disjoint sums compatible with fiber products if there exists an ordinal $\beta$ such that the sieves corresponding to the covering families of size strictly smaller than $\beta$ generate the topology (i.e. they are cofinal among every sieve), and if disjoint sums of less than $\beta$ elements exist in $\X$ and are compatible with fiber products (i.e. we have a distributive property for the fiber product of two disjoint sums over an object of $\X$). \\
For example, if $\X$ is quasi-compact we can take $\beta=\omega$ and the condition is then that there are finite disjoint sums. \\
We will say that the disjoint sums are covering if for every family $\mathfrak{U}$ closed under disjoint sums, the family $\mathfrak{U}$ is a family covering $X$.
\end{defn}
Let $\X$ be a site equipped with a sheaf of rings $\Ox$. We suppose that $\X$ admits fiber products and that it admits enough disjoint sums compatible with fiber products. We write for $X\in \X$, $\cc[X]_{\Ox}$ the category of complexes of sheaves of $\Ox$-modules over $\X/X$. \\
Let $\mr{qis}(X)\subset \cc[X]_{\Ox}$ denote the class of quasi-isomorphisms, and we write $\cc[X]_{\Ox}^{[0,\infty)}$ for the subcategory of those complexes which are bounded below and we consider a cofibrantly generated model category structure on $\cc[X]_{\Ox}^{[0,\infty)}$ (for example, the so-called injective model category structure on chain complexes). We write $L(\cc[X]_{\Ox}^{[0,\infty)},\mr{qis}(X))$ for the simplicial category obtained by simplicial localization at the class of morphisms $\mr{qis}(X)$.\\
We can now write down the main corollary of Hirschowitz-Simpson:
\begin{cor}\label{cor:descent}\cite[Corollary 21.1]{hirschowitz1998descente}
Suppose $\X$ admits fiber products and enough compatible disjoint sums. Then the presheaf of simplicial categories 
\[ L(\mathscr{C}h_{\Ox}^{[0,\infty)}):=(X\mapsto L(\cc[X]_{\Ox}^{[0,\infty)},\mr{qis}(X))) \]
is a Segal 1-stack over $X$ which we call the Segal 1-stack of modules of complexes over $(\X,\Ox)$.
\end{cor}
We refer the reader to \cite[Chapter 3, Theorem 1.1]{hirschowitz1998descente} and \cite{simpson2011homotopy} for a detailed account of Segal n-stacks and Segal n-categories and their homotopy theory. 
\begin{obs}\label{remark:inftymodels}
The language of higher Segal categories and higher Segal stacks of Hirschowitz and Simpson is, while of importance from the point of view of the foundations of higher category theory, not the standard modern framework in which higher categories and higher stacks is often presented in the literature. It is however understood by the experts that this language is for the most part equivalent to constructions performed in terms of quasi-categories. \\
In \cite{toen2007moduli} Toën and Vaquié construct a derived and more general version of this Segal 1-stack with a focus in the language of model categories and dg-categories. See \cite[Definition 3.28]{toen2007moduli} and their following remarks.  We can also point out to \cite{meadows2020descent} for a modern treatment of Corollary \ref{cor:descent}.\\
\end{obs}
The original motivation of Hirschowitz and Simpson of treating the problem of gluing of (bounded) chain complexes up to quasi-isomorphism is included in this theorem. When specializing to $\X=\mr{Sch}$ the Zariski site of schemes this result implies that at the level of dg-categories we can carry out descent. Very concretely and with the comments in Remark \ref{remark:inftymodels} in mind, we the following statement which justifies the assumption $S=\mr{Spec} B$ all through this work.
\begin{cor}\label{cor:affineimplies}
    Let $X$ be a quasi-compact and quasi-separated scheme, the $\infty$-stack of bounded above derived categories over $X$ in Corollary  \ref{cor:descent} exists and is defined by finite affine covers.
\end{cor}
As the key step in the proof of our main Theorem revolves around the behavior around a point $s\in S$ we will now constrain ourselves to that particular situation. For self-consistency and to emphasize the use of dg-categories we switch the nomenclature in Corollary \ref{cor:descent} when dealing with this concrete geometric context: 
\begin{defn}\label{def:bape}
    Let $X\to S$ be a smooth projective faithfully flat variety over $S$. We write $\bape$ for the $\infty$-stack of dg-categories defined by the assignment affine open subsets
     \[ U \mapsto \mr{Perf}^{dg}(X_{U}) \]
\end{defn}
\begin{obs}\label{obs:stackiffication}
We incur in a slight abuse of notation, but the concerned reader should understand terms like "$\infty$-stack of dg-categories" as functors with target in a $(\infty,1)$-category or an appropriate model category of simplicial sets.\\
We will frequently use the fact that one can define a sheaf $\F$ of dg-categories, objects in dg-categories, and $\infty$-groupoids by specifying the value of the presheaf $F$ at affine opens since these form a basis of open subsets in the usual sense. \\
With these values defined, we can obtain the associated sheaf by sheaffification
\[ \F(V):= \underset{U^{\bullet}\to V}{\mr{lim}} F(U^{\bullet})\]
Here the limit is taken in an appropriate sense. When working with sheaves of $(\infty,1)$-categories, see \cite[Section 6.5.3]{lurie2009higher} or alternatively when working with sheaves of simplicial sets see \cite[Definition 2.1.2]{toen2010simplicial}.
\end{obs}
We can for example define:
\begin{defn}\label{def:sheaftensor}
    Let $\bape$ and $\bbpe$ as in the statement of Corollary \ref{thm:descentformodcats} then we define the tensor product sheaf defined as the associated presheaf which on affine open subsets $U\subset S$ we describe as
    \[ U\mapsto \bape(U)\tofdgc \bbpe(U) \]
    and we denote it by $\bape\tofsh\bbpe$.
\end{defn}
The rest of this section seeks to justify the existence of the $\infty$-stack of tensor bifunctor structures on $X\to S$. This is simply a special case of the stack $\mathbb{R}Perf$ of \cite{toen2007moduli} and which one should think of as a derived version of the 1-stack $\underline{Coh}$ defined by Laumon-Moret Bailly in \cite{lmb}. 
\begin{cor}\label{thm:descentformodcats}
Let $n,m\in \mathbb{N}$, $X\to S$, and $\bape$ as above. Write $\bbpe$ for the stack of dg-categories on a smooth projective flat variety $Y\to S$ as described in Definition \ref{def:bape}. The presheaf of $\infty$-categories 
\[ \ihom{\bape^{\tofdgc n}}{\bbpe^{\tofdgc m}}:=(U\mapsto \mr{Int}(\bape(U)^{\tofdgc n}\tofdgc (\bbpe(U))^{op\,\tofdgc m}-\mr{Mod}^{\mr{rqr}})) \]
induces a stack over S.
\begin{proof}
   The result follows from Corollary \ref{cor:descent}, the construction of Definition \ref{def:sheaftensor}, and the Morita theorem for dg-categories characterizing the internal Hom of the $\infty$-category dgCat as modules over $\bape^{\tofdgc n}\tofdgc \bbpe^{op\,\tofdgc m}$. \\
   Let us sketch the construction along two open affine subsets $U,V\subseteq S$ with $U\cup V=S$. Corollary \ref{cor:descent} implies there is a cartesian square: 
   \[ \xymatrix{\bape(S)\ar[r] \ar[d] & \bape(U)\ar[d] \\ \bape(V) \ar[r] & \bape(U\cap V)}\]
   Which we can tensor by $\bape^{\tofsh n-1}$ and $\bbpe^{op\,\tofsh m}$ to obtain a square:
   \[ \xymatrix{\bape^{\tofsh n}(S)\tofsh \bbpe^{op\,\tofsh m}(S) \ar[r] \ar[d] & \bape(U)^{\tofdgc n}\tofdgc \bbpe(U)^{op\,\tofdgc m} \ar[d] \\ \bape(V)^{\tofdgc n}\tofdgc \bbpe(V)^{op\,\tofdgc m} \ar[r]  & \bape(U\cap V)^{\tofdgc n}\tofdgc \bbpe(U\cap V)^{op\,\tofdgc m}}.\]
   By applying the functor $\mr{Int}(\_)-\mr{Mod}^{rqr}$ we get a diagram
   \[ \xymatrix{\ihom{\bape(S)^{\tofdgc n}}{ \bbpe(S)^{\tofdgc m}}  & \ihom{\bape(U)^{\tofdgc n}}{\bbpe(U)^{\tofdgc m}}  \ar[l] \\ \ihom{\bape(V)^{\tofdgc n}}{\bbpe(V)^{\tofdgc m}} \ar[u] & \ihom{\bape(U\cap V)^{\tofdgc n}}{\bbpe(U\cap V)^{\tofdgc m}} \ar[u] \ar[l] }\]
   Which by the Morita theorem corresponds to the cartesian square:
   \[ \xymatrix{\mr{Perf}^{dg}(X\times_{S} Y) \ar[r] \ar[d] & \mr{Perf}^{dg}(X_{U}\times_{S} Y_{U}) \ar[d] \\ \mr{Perf}^{dg}(X_{V}\times_{S} Y_{V})\ar[r] & \mr{Perf}^{dg}((X_{U}\cap X_{V})\times_{S} (Y_{U}\times Y_{V})}\]
   As pointed out by Toën and Vaquié in \cite[Definition 3.28]{toen2007moduli} one can think of the stack $\ihom{\bape^{\tofdgc n}}{\bbpe^{\tofdgc m}}$ (or $\mathbb{R}\underline{Perf}(X^{\times n}\times Y^{\times m})$ in the notation of Toën-Vaquié) as the internal Hom stack
   \[ Map(X^{\times n}\times_{S} Y^{\times m}, \bape\tofsh \bbpe).\]
\end{proof}
\end{cor}
In particular by setting $X=Y$ and $n=2, m=1$, we obtain the existence of a stack of modules over $\bape$ which parameterize 2-fold dg-bimodules. We will make use of this construction in our next section.

\section{Stacks of tensor product structures}\label{sec:inftystack}
We've come to the main part of our current work. Here we will take advantage of the constructions from the previous section and consider a $\infty$-stack of modules which we will think of as a stack of tensor structures over $X\to S=\mr{Spec} B$. \\
We are concretely interested in understanding the passing from the behaviour at a fiber $X_{s}$ along a point $s\in S$ of this stack to the behaviour around an affine neighborhood. \\
To recall our initial set-up, we have a TTC $(\boxtimes,\mathbb{1})$ on the derived category $\der[X\to S]$, the bifunctor associated to $\boxtimes$ has a dg-enhancement in the form of a (right quasi-representable) 2-fold dg-bifunctor $\Gamma\in \mr{Bimod}^{2}_{dg}(\bape(X))$. In Section \ref{sec:xtofibers} we added the bifunctor $\boxtimes$ was reasonable (Definition \ref{def:reasonable}) with respect to the adjoint pair $I_{s}^{\ast}\dashv I_{s}\,_{\ast}$ for a particular fixed point $s\in S$ in order to be able to discuss an induced bifunctor on the derived category $\der[X_{s}=X\times_{S} \{s\}]$ of the fiber $X_{s}$. \\
Now we want to go in the other direction, our first step is then to define a sheaf on $X\to S$ induced by the global object $\Gamma$. \\
\begin{defn}\label{def:georeasonable}
We say that a MTC $(\boxtimes, \mathbb{1})$ is geometrically reasonable if for every open affine neighborhood $U$ of $S$ with $I_{U}:X_{U}\hookrightarrow X$, it is reasonable with respect to $I_{U}^{\ast}\dashv I_{U}\,_{\ast}$, and additionally if for every $i:U\hookrightarrow V$ the induced restriction dg-functors $\Phi_{VU}:\bape(V)\to \bape(U)$ are monoidal in the sense that we have coherent equivalences
\[ \Phi_{VU}(\_\boxtimes_{V}\_)\simeq (\_)\boxtimes_{U}(\_) \]
Where $\boxtimes_{V}:=\boxtimes_{I_{V}^{\ast}I_{V}\,_{\ast}}$ (and where $\boxtimes_{U}$ is analogously defined).
\end{defn}
Since every pullback and pushforward derived functor induced by inclusions is of Fourier-Mukai type (see \cite[Example 5.4]{huybrechts2006fourier}) then they all have dg-enhancements in the form of dg-bimodules and so do their compositions with other functors of Fourier-Mukai type. \\
We will write $\Gamma_{U}$ for the dg-enhancements functors $\boxtimes_{U}$ at open affine subsets, so that we have families $\Gamma_{U}\in \mr{Bimod}_{dg}^{2}(\bape(U))$ indexed by $U\in S$. \\
We want to construct a sheaf $\bGamma$ of $\infty$-groupoids induced locally by these 2-fold dg-bimodules $\Gamma_{U}$.
\begin{prop}\label{prop:gammaisheaf}
The sheaf described above exists and is determined by the values of the presheaf at affine open subsets $U\subset S$ under the condition that $\boxtimes$ is geometrically reasonable. 
    \begin{proof}
    We will unwrap the description above. \\
    If $U\subseteq S$ is an affine open we will write $\bGamma(U)$ for the $\infty$-groupoid
    \[ \{\Gamma_{U}\in \mr{Bimod}_{dg}^{2}(\bape(U))\mid \Gamma_{U}\,\,\mr{enhances}\,\boxtimes_{U}\}\]
    We take this $\infty$-groupoid as the equivalence classes in the $\infty$-category associated to the model category structure on $\mr{Bimod}_{dg}^{2}(\bape(U))$. \\
    It is clear that there are compatible actions by $\bape(U)$ on the left and $\bape(U)^{op}$ on the right of $\bape(U)$ by which we mean that there are morphisms in the $\infty$-category dgCat
    \[ \bGamma(U)\otimes \bape(U)\tofdgc \bape(U) \to \bape(U)^{op}\otimes\cc[k] \]
    Here the tensor product $\bape(U)\otimes (\bape(U)\tofdgc \bape(U))$ denotes the tensoring action of the presentable $\infty$-category $\bGamma(U)$ with the dg-category $(\bape(U)\tofdgc\bape(U))$. \\
    For any $i:U\hookrightarrow V\subseteq S$ the restriction $\phi_{VU}$ induces a restriction
    \[\phi_{VU}:\bGamma(V)\to \bGamma(U) \]
    We obtain the following commutative diagrams
    \begin{align*}
    \xymatrix{\bGamma(V)\otimes\bape(V)\tofdgc\bape(V) \ar[d]^{\phi_{VU}\otimes\Phi_{VU}\otimes\Phi_{VU}} \ar[r]^{m_{V}}  & \bape(V)^{op}\otimes \cc[k] \ar[d]^{\Phi_{VU}\otimes Id} \\ \bGamma(U)\otimes\bape(U)\tofdgc\bape(U) \ar[r]^{m_{U}} & \bape(U)^{op}\otimes\cc[k]}
    \end{align*}
The commutativity of this diagram is satisfied by the assumption that $\boxtimes$ was geometrically reasonable. \\
By realizing the 2-fold dg-bimodules $\Gamma_{U}$ as objects in $\mathbb{R}\underline{Hom}(\bape(U)\tofdgc \bape(U),\bape(U))$ then by the construction of Corollary \ref{thm:descentformodcats}, we know that these objects glue into a $\infty$-stack on $\infty$-groupoids of dg-modules over $\bape$. 
\end{proof}
\end{prop}
As with $\bape$ we will again abuse our notation and keep writing $\bGamma$ for the sheaf defined above.  
\begin{obs}
We would also like to mention the ideas of six functor formalisms developed by Clausen-Mann-Scholze in \cite{scholze, mann2022p} as it seems our constructions here can be seen as very concrete cases of their theory. In particular we are focusing on noncommutative derived schemes and exploiting the derived Morita correspondence of Toën for our specific applications. In their framework and their language, our construction here corresponds to considering 6-functor formalisms 
\[ D_{\boxtimes}:\mr{Corr(Sch/X,E)}\to \mr{Mon(Cat_{\infty})}\]
For the reader's convenience we quickly recall Mann's definition borrowed from \cite{scholze}.\\
Let $C$ be a $\infty$-category with finite limits and let $E$ a class of morphisms stable under pullbacks and compositions. The category of correspondences $Corr(C,E)$ has the same objects as $C$, morphism objects and the composition $Corr(C,E)\times Corr(C,E)\to Corr(C,E)$ are defined by $\infty$-groupoids of spans. The category $Corr(C,E)$ has a monoidal structure given by the cartesian monoidal structure on $C$. \\
\begin{defn}\label{def:3ff}
A 3-functor formalism is a lax symmetric monoidal functor $D:Corr(C,E)\to Cat_{\infty}$.
\end{defn}
Here $Cat_{\infty}$ denotes the $\infty$-category of $\infty$-categories. As shown by Mann, this definition encodes the assignment to each object of $C$, a symmetric monoidal $\infty$-category $D(C)$, and the existence of pullback and direct image functors $f^{\ast}:D(X)\to D(Y)$, $f_{!}:D(Y)\to D(X)$ for every morphism $f:Y\to X$.
\begin{defn}\label{def:6ff}
    A 6-functor formalism is a 3-functor formalism such that the monoidal structures on each $D(X)$ is closed, and the pullback and direct image morphisms all have right adjoints.
\end{defn}
What these additional conditions guarantee is the existence of a base change and projection formulas. \\
In this language, our setting corresponds to considering the relative big Zariski affine site $Aff/X$, and so assignments $U\mapsto D(U)$ where the underlying $\infty$-category is the dg-category $\bape(U)$ and the monoidal structure is one enhancing the $\boxtimes_{U}$ where $\boxtimes$ is a geometrically reasonable MTC structure (Definition \ref{def:reasonable}. \\
From this perspective one immediate success of Mann's definition is that it encodes some of the conditions we require our tensor triangulated categories to have in a convenient and compact definition. It would be of interest to come back to this perspective and rephrase our work here in this language and generality. \\
Some of the required conditions in this formalism are automatically verified in our case by the particular context of consideration. Other conditions like a version of the projection formula comes from the geometric reasonableness we ask our tensor structure to have. 
\end{obs}
\begin{defn}\label{defn:inftystalk}
Let $\Ox(Y)$ be the category of open subsets of a topological space $Y$. Consider a presheaf $\F:\Ox(Y)^{op}\to \C$ with values in some presentable $\infty$-category $\C$.\\
Let $y\in Y$, the stalk of $\F$ at $y$ is the object $\F_{y}$ of $\C$ defined as
\[ \underset{y\in U\subseteq Y}{\mr{Colim}} \F(U) \]
\end{defn}
This is essentially the usual definition of a stalk of sheaves except for the fact that we are considering a (filtered) colimit in the $\infty$-categorical sense. \\
\begin{lemma}\label{lemma:tensorstalk}
Let $X\to S$, and $\bape$ as above. Consider a point $s\in S$, then for any $n\in \mathbb{N}$ we have
\[(\bape^{\tofsh n})_{s}\simeq ((\bape)_{s})^{\tofdgc n}\]
\begin{proof}
This follows from the fact that the tensor product of dg-categories is closed and symmetric and thus commutes with colimits. By the definition of the tensor product of sheaves of dg-categories we just need to refine the directed system $s\in U\subseteq S$ to consist only in affine opens $U\ni s$.
\end{proof}
\end{lemma}
\begin{lemma}\label{lemma:opstalk}
Let $X\to S$ and $\bape$ be as above. Then $(\bape)_{s}^{op}\simeq (\bape^{op})_{s}$.
\begin{proof}
By simple computation, 
\[ (\bape)_{s}^{op}:=(\underset{s\in U\subseteq S}{\mr{Colim}}\bape(U))^{op} \simeq \underset{s\in U\subseteq S}{\mr{Colim}}\bape(U)^{op}=:(\bape^{op})_{s}\]
\end{proof}
\end{lemma}
\begin{lemma}
\label{lemma:bGammastalk}
Let $X\to S$, $\bape$ and $\bGamma$ as above. Take $s\in S$, then $\bGamma_{s}$ is a 2-fold dg-bimodule over $(\bape)_{s}$. 
\begin{proof}
By definition
\[ \bGamma_{s}:=\underset{s\in U\subseteq S}{\mr{Colim}}\bGamma(U)\]
We have that each $\bGamma(U)$ is a $\infty$-groupoid of 2-fold dg-bimodule over $\bape(U)$. The result follows from noticing that the action compatibility
\begin{align*}
\xymatrix{ \bGamma(V)\otimes\bape(V)\tofdgc\bape(V)\ar[r] \ar[d] & \bape(V)^{op}\otimes \cc[k] \ar[d] \\ \bGamma(U)\otimes\bape(U)\tofdgc\bape(U)\ar[r] & \bape(U)^{op}\otimes\cc[k]}
\end{align*}
commutes with colimits by Lemmas \ref{lemma:tensorstalk} and \ref{lemma:opstalk} and the fact that the tensoring of the presentable $\infty$-categories $\bGamma(U)$ with a dg-category commutes with colimits (again using that this monoidal structure is closed and symmetric). 
\end{proof}
\end{lemma}
The next immediate goal of this section is to write a precise relationship between the stalk dg-category $(\bape)_{s}$ around a point $s\in S$ and the dg-category of perfect complexes of the fiber $X_{s}:=X\times_{S}\{s\}$. \\
Furthermore, for $i:X_{s}\to X$ the fiber inclusion and $\boxtimes$ a $S$-linear TTC structure on $\mr{Perf}(X)$ reasonable with respect to $i^{\ast}\dashv i_{\ast}$ we would like to exhibit $\boxtimes_{i^{\ast}i_{\ast}}$ as canonically induced by $\bGamma_{s}$. \\
\\
Before moving on to the main results of this section, let us briefly recall some important constructions in the general theory of dg-categories and some properties of these that we will be using for one of our main results. \\
In Definition \ref{def:dgquotient} we described the dg-quotient of a dg-category $\A$ by a full dg-sub-category $\B\subset \A$. Let us expand more on the details of this construction. \\
In this situation where we are working over a field, the Hom chain complex objects $Hom_{\A/\B}(X,Y)$ for any two objects $X,Y\in \A/\B$ can be described as a span
\begin{align}\label{defn:quotientspan}
    \xymatrix{X=X_{0} \ar[dr]^{f_{0}} & & X_{2} \ar[dr]^{f_{2}} \ar[ld]^{\epsilon} & & X_{4} \ar[dl]^{\epsilon} \\ & X_{1} & & X_{3}  } \dots \xymatrix{X_{n-2} \ar[dr]^{f_{n-2}} & & Y=X_{n} \ar[ld]^{\epsilon} \\ & X_{n-1} & }
\end{align}
So that the graded $k$-module associated to $Hom_{\A/\B}(X,Y)$ is, in degree $n$, a product $Hom_{\A}(X=X_{0},X_{1})\otimes k[1]\otimes Hom_{\A}(X_{1},X_{2})\otimes k[1]\otimes \dots \otimes Hom_{\A}(X_{n-1},X_{n}=Y)$, for some objects $X_{i}\in \B$ for $0 < i < n$. \\ 
One can use this to find an explicit construction of the Hom chain complexes by recalling that $\epsilon$ is such that $d(\epsilon)=1$ for any object in $\B$. See \cite[Construction 3.1]{drinfeld2004dg} for the concrete details of this description. \\
Tabuada later on gives the following useful universal property of the quotient in the homotopy category $H_{qe}$. \\
First, recall that one says a dg-functor $\F:\A\to \M$ is said to annihilate the dg-subcategory $\B\subseteq \A$ if for any $X\in \B$, the identity morphism $1_{X}:X\to X$ vanishes in $H^{0}(\M)$ by $H^{0}(\F)$. \\
With this definition we can phrase Tabuada's theorem.
\begin{thm}(\cite[Theorem 4.0.1]{tabuada2010drinfeld})\label{thm:tabuadaupquotient}
Let $\A$ be a dg-category and $\B\subseteq \A$ be a full dg-subcategory of $\A$. For every dg-category $\M$, the Drinfeld dg-quotient $Q:\A\to \A/\B$ induces a bijection 
\[ H_{qe}(\A/\B,\M) \overset{\simeq}{\longrightarrow} H_{qe\B}(\A,\M) \]
Where $H_{qe\B}(\A,\M)$ denotes the set of morphisms which annihilate $\B$.
\end{thm}
As pointed out by Tabuada, this characterizes the Drinfeld dg-quotient as a homotopy cofiber sequence in $H_{qe}$. \\
The next theorem is a reframing to our dg-category language of Lemma 3.3(a) of \cite{thomason1997classification}. There, Thomason establishes a relationship between acyclic complexes $E^{\bullet}$ of sheaves of $\Ox_{X,x}$-modules at the stalk of a point $x\in X$ where $X$ is a quasi-compact quasi-separated scheme, and the acyclic complexes $E^{\bullet}\dteee k(x)$ by the concept of homological support. \\
For convenience we include Thomason's lemma, as originally written:
\begin{lemma}\cite[Lemma 3.3(a)]{thomason1997classification}\label{lemma:thomasonsuppo}
Let $X$ be a quasi-compact and quasi-separated scheme. Let $E^{\bullet}$ be a perfect complex on $X$. Then for any $x\in X$, $E_{x}$ is an acyclic complex of $\Ox_{X,x}$-modules if and only if $E\dteee k(x)$ is an acyclic complex of $k(x)$-modules.
\end{lemma}
Now consider the inclusion 
\[ FI_{U}:X\times_{S}\{s\}=X_{s}\hookrightarrow X_{U}:=X\times_{S}U\]
And the induced derived pullback
\[ FI_{U}^{\ast}:\bape(U)\to \bape(s)\]
Let us denote by $HC_{U}:\bape(U)\to (\bape)_{s}$ the functor given by the colimit construction. By the universal property of this colimit, we know that there is a map 
\begin{align}\label{eq:UP}
UP:(\bape)_{s}\to \bape(s)
\end{align}
such that $UP\circ HC\simeq FI_{U}^{\ast}$. \\
\begin{lemma}\label{lemma:thomasonsupp}
    Let $X\to S$, and $\bape$ be as above. Let $s\in S$, then $E^{\bullet}_{s}\in (\bape)_{s}$ is acyclic if and only if $UP(E^{\bullet})\in \bape(s)$ is acyclic.
    \begin{proof}
        The proof follows the classical reasoning present in Thomason's original proof. \\
        Let us consider the spectral sequence
        \begin{align}\label{eq:sps}
            E_{2}^{p+q}:=H^{p}(UP(H^{q}(E^{\bullet}_{s}))) \Rightarrow H^{p+q}UP(E_{s}^{\bullet})
        \end{align}
        So as usual, if $E_{s}^{\bullet}$ is acyclic, so is $UP(E^{\bullet}_{s})$. \\ 
        Conversely, suppose $E_{s}^{\bullet}$ is not acyclic. As we deal with perfect complexes, $E_{s}^{\bullet}$ is cohomologically bounded and there is a least $N$ with $H^{N}(E^{\bullet}_{s})\not=0$. We look at the quadrant of $E_{2}^{p,q}$ where $p>0$ and $q> N-1$, out of which we know every object is 0.
        \begin{align*}
            \xymatrix @C=1pc @R=1pc{ & \vdots & \vdots & \vdots &  \\ 0 &  H^{0}(UP(H^{q+2}E^{\bullet}_{s})) \ar[drr] & H^{1}(UP(H^{q+2}E^{\bullet}_{s})) \ar[drr] & H^{2}(UP(H^{q+2}E^{\bullet}_{s})) & \dots \\ 0 & H^{0}(UP(H^{q+1}E^{\bullet}_{s})) \ar[drr]  & H^{1}(UP(H^{q+1}E^{\bullet}_{s})) \ar[drr] & H^{2}(UP(H^{q+1}E^{\bullet}_{s}))  & \dots \\0 &  H^{0}(UP(H^{q}E^{\bullet}_{s})) \ar[drr] & H^{1}(UP(H^{q}E^{\bullet}_{s})) \ar[drr]  & H^{2}(UP(H^{q}E^{\bullet}_{s})) &\dots \\  & 0 & 0 & 0 & \dots}
        \end{align*}
        Therefore, there is an isomorphism $H^{0}(UP(H^{N}E_{s}^{op}))\simeq H^{N}E_{s}^{\bullet}$. \\
        Furthermore, $H^{N}E^{\bullet}_{s}$ is an object in $(\bape)_{s}$ and so we can write as a colimit $\underset{s\in U}{\mr{Colim}}H^{N}E^{\bullet}_{U}$ for objects $E_{U}^{\bullet}\in \bape(U)$. Again, the objects $H^{N}E^{\bullet}_{U}$ belong to $\bape(U)$ too and so to conclude the proof, Corollary \ref{cor:nakayama} of the dg-Nakayama Lemma \ref{lemma:nakayama} below applied to  $H^{N}E^{\bullet}_{U}$ implies $UP(H^{N}E_{s}^{\bullet})\not\simeq 0$ and so $UP(E_{s}^{\bullet})$ is not acyclic.
    \end{proof}
\end{lemma}
In \cite{zimmermann2023differential}, Zimmermann gave a proof of a differential graded version of the classical Nakayama Lemma. We need the following definition to phrase the concrete statement of the theorem. 
\begin{defn}\label{def:dgrad2}
Let $A$ be a dg-algebra, the differentially graded two-sided radical $\mr{dgrad}_{2}(A)$ is the intersection of all the annihilators of dg-simple differentially graded left $A$-modules. In other words
\[ \mr{dgrad}_{2}(A):=\bigcap_{S}\mr{Ann}(S):=\bigcap_{S}\{r\in A\mid rS=0\}\]
\end{defn}
Where dg-simple means that the module has no non-trivial dg-submodules. \\
\begin{lemma}\cite[Lemma 4.27]{zimmermann2023differential}\label{lemma:nakayama}
Let $A$ be a dg-algebra. Let $M, N$ be dg-A-modules such that $N\subseteq M$. Assume that $M$ is finitely generated as a dg-module. Then $\mr{dgrad}_{2}(A)M=M$ implies $M=0$, and $N+\mr{dgrad}_{2}(A)M=M$ implies $M=N$.
\end{lemma}
\begin{obs}\label{obs:orlovjacobson}
    An alternative Nakayama Lemma was proven by Goodbody in \cite[Theorem 3.3]{goodbody2024reflecting} for fd dg-algebras using a version of the Jacobson radical introduced first by Orlov. It seems that the language used by Goodbody and Orlov fits better with the homotopical nature of our work, but unfortunately and as pointed out by Goodbody \cite[Remark 2.4]{goodbody2024reflecting}, it is not the case that the dg-algebras coming from algebraic varieties are necessarily finite dimensional.
\end{obs}
Zimmermann's Nakayama's Lemma has the usual corollaries for formal reasons. In particular the usual geometric form of the lemma:
\begin{cor}\label{cor:nakayama}
    Let $A$ be a dg-algebra and $M$ a finitely generated dg-A-module. Let $\{m_{1},\dots,m_{k}\}\subset M$ and suppose $M/\mr{dgrad}_{2}(A)M$ is generated as a $A$-module by the images of $\{m_{i}\}$  under the quotient morphism $M\to M/\mr{dgrad}_{2}(A)M$, then $\{m_{i}\}$ generate $M$ as a $A$-module. 
\end{cor}
        
 Let us recall,
\begin{defn}\label{def:hsupport}
    Let $F^{\bullet}$ be a complex of sheaves on some space $X$. The homological support $\mr{Supph}F^{\bullet}$ of $F^{\bullet}$ is defined as $\bigcup_{i\in \mathbb{Z}} \mr{Supp}\mathcal{H}^{i}F^{\bullet}$.
\end{defn}
\begin{defn}\label{def:hsupp}
The support dg-subcategory $\mathfrak{m}_{U}$ of $\bape(U)$ is the full dg-category spanned by 
\[\{P^{\bullet}\in \bape(U)\mid \mr{Supph}P^{\bullet}=\emptyset\}\]
\end{defn}
Notice that since $\mr{Supph}$ is a cohomological invariant then $\mathfrak{m}_{U}$ is well defined.  We write $\mathfrak{m}_{s}$ for $\underset{s\in U}{\mr{\mr{Colim}}}\;\mathfrak{m}_{U}$. \\
We now assemble everything to establish a relationship between the dg-categry of perfect complexes on the fiber at a point $s\in S$, and $(\bape)_{s}$, mimicking the classical situation in algebraic geometry. 
\begin{thm}\label{thm:fiberisq}
Let $X\to S$, and $\bape$ be as above. Let $s\in S$, then there is an equivalence of dg-categories
\[ \bape(s)\simeq (\bape)_{s}/\mathfrak{m}_{s}. \]
\begin{proof}
As each dg-category $\mathfrak{m}_{U}$ is a dg-subcategory of $\bape(U)$ for every $U\ni s$ then it is clear that $\mathfrak{m}_{s}$ is a full dg-subcategory of $(\bape)_{s}$ and the quotient is well defined. \\
Assembling everything together we have the following diagrams for each $U\subset S$
\begin{align*}
\xymatrix{ &\bape(U) \ar[dl]^{HC_{U}} \ar[dr]^{FI_{U}^{\ast}} & \\ (\bape)_{s}\ar[rr]^{UP} \ar[dr]^{q} & & \bape(s) \\ & (\bape)_{s}/\mathfrak{m}_{s} & }
\end{align*}

Where the upper triangle commutes. \\
We want to show that $FI_{U}^{\ast}$ annihilates $\mathfrak{m}_{U}$ for every $U$. But this is true, by Lemma \ref{lemma:thomasonsupp} we have that $P^{\bullet}\in (\bape)_{s}$ is acyclic if and only if $UP(P^{\bullet})\in \bape(s)$ is acyclic. We can then refine the set of open subsets $U\subseteq S$ such that our claim holds and this has no effect on the calculation of the colimits. \\
Since $FI_{U}^{\ast}$ annihilates $\mathfrak{m}_{U}$ then we have maps $DQ_{U}:\bape(U)/\mathfrak{m}_{U}\to \bape(s)$ such that $DQ_{U}\circ q_{U}\simeq FI_{U}^{\ast}$, where $q_{U}$ denotes the quotient map. And therefore by the colimit universal property, of $\mathfrak{m}_{s}$, we obtain a morphism 
\[ TE:(\bape)_{s}/\mathfrak{m}_{s} \to \bape(s) \]
Additionally, by the universal property of the Drinfeld dg-quotient (Theorem \ref{thm:tabuadaupquotient}) we obtain a map $(\bape)_{s}/\mathfrak{m}_{s}\to \bape(s)$ which coincides with $TE$. \\
This can be seen from the commutativity of the following diagram
\begin{align}\label{diagram:allcommutes}
\xymatrix{ &\bape(U) \ar[dl]_{HC_{U}} \ar[d]^{FI_{U}^{\ast}} \ar[dr]^{q_{U}} & \\ (\bape)_{s}\ar[r]^{UP} \ar[dr]_{q} & \bape(s) & \bape(U)/\mathfrak{m}_{U} \ar[l]^{DQ_{U}} \ar[dl]^{SHC_{U}} \\ & (\bape)_{s}/\mathfrak{m}_{s} \ar[u]^{TE} & }
\end{align}
That $TE$ is an equivalence follows again from Lemma \ref{lemma:thomasonsupp} since the Drinfeld dg-quotient is a cofiber sequence.
\end{proof}
\end{thm}
\begin{obs}\label{obs:bapes}
The reader might notice a slight abuse of notation in the term $\bape(s):=\bape(X_{s})$ as $\bape$ is only defined on open subsets $U\subseteq S$ (alternatively on varieties $X_{U}:=X\times_{S}U$). What we mean then by $\bape(s)$ is simply a dg-enhancement of the category of perfect complexes $\mr{Perf}(X_{s})$. No step in the considerations above depend on whether $\bape(s)$ is a well defined term for the stack $\bape$ and so we believe keeping this notation to mirror the usual algebraic-geometric one is better for our exposition despite the discrepancy in the notation $(\bape)_{s}\not\simeq \mr{Perf}^{dg}(X_{s})\simeq \bape(s)$.
\end{obs}
Next we show that one can argue about $\bGamma$ in the same way we argued about $\bape$ above.  
\begin{lemma}\label{cor:bgammasisq}
The dg subcategory $\ihom{\bape(s)^{\tofdgc 2}}{\bape(s)}_{ac}$ given by the essential image of the map $\ihom{\bape(s)^{\tofdgc 2}}{(\bape)_{s}}\to \ihom{\bape(s)^{\tofdgc n}}{\bape(s)}$ is a Drinfield dg-quotient of $\ihom{\bape(s)^{\tofdgc n}}{(\bape)_{s}}$
\begin{proof}
    We apply $\ihom{\bape(s)^{\tofdgc 2}}{\_}$ to Diagram \ref{diagram:allcommutes} and obtain
    \[ \xymatrix{ & \ihom{\bape(s)^{\tofdgc 2}}{\bape(U)} \ar[d] \ar[dl] \\ \ihom{\bape(s)^{\tofdgc 2}}{(\bape)_{s}}  \ar[r]  & \ihom{\bape(s)^{\tofdgc 2}}{\bape(s)}} \]
    With the maps being induced by composition with the respective map of Diagram \ref{diagram:allcommutes}. However, note that since the functor $\ihom{\bape(s)^{\tofdgc 2}}{\_}$ commutes with filtered colimits, the morphism $\ihom{\bape(s)^{\tofdgc 2}}{(\bape)_{s}}\to \ihom{\bape(s)^{\tofdgc 2}}{\bape(s)}$ coincides with the map induced by the universal product of this filtered colimit. \\
    Let us write $H\mathfrak{m}_{U}$ for the full dg subcategory spanned by
    \[\{F\mid F(P^{\bullet},Q^{\bullet})\in \mathfrak{m}\;\forall\;P^{\bullet},Q^{\bullet}\in \bape(s) \}. \]
    Notice then that the morphism $\ihom{\bape(s)^{\tofdgc 2}}{\bape(U)}\to \ihom{\bape(s)^{\tofdgc 2}}{\bape(s)}$ annihilates $H\mathfrak{m}_{U}$. \\
    Indeed, if $F\in H\mathfrak{m}_{U}$ then its image in $\ihom{\bape(s)^{\tofdgc}}{\bape(s)}$ must be acyclic for every $P^{\bullet},Q^{\bullet}\in \bape(s)$. we obtain the following diagram 
    \begin{align}\label{diagram:ihomallcommutes}
  \xymatrix@C=1mm{ & \ihom{\bape(s)^{\tofdgc 2}}{\bape(U)} \ar[d] \ar[dl] 
 \ar[dr]^{Hq_{U}} & \\ \ihom{\bape(s)^{\tofdgc 2}}{(\bape)_{s}}  \ar[r] \ar[dr]  & \ihom{\bape(s)^{\tofdgc 2}}{\bape(s)}_{ac} & \ihom{\bape(s)^{\tofdgc 2}}{\bape(U)}/H\mathfrak{m}_{U} \ar[l]^{HDQ_{U}} \ar[ld] \\ & \ihom{\bape(s)^{\tofdgc 2}}{(\bape)_{s}}/H\mathfrak{m}_{s} \ar[u]^{TE} &}
    \end{align}
    Let us explain, the map $Hq_{U}$ is the quotient map, $HDQ_{U}$ is the induced map by the universal property of the Drinfeld quotient, and the map \[TE:\ihom{\bape(s)^{\tofdgc 2}}{(\bape)_{s}}/H\mathfrak{m}_{s}\to \ihom{\bape(s)^{\tofdgc 2}}{\ihom(s)}\] 
    is the equivalence we are looking for, where $H\mathfrak{m}_{s}$ is the colimit of the dg subcategories $H\mathfrak{m}_{U}$.  
\end{proof}
\end{lemma}
We can now phrase the main theorem of this work:
\begin{thm}\label{thm:mainthm22222222}
    Let $X\to S$ be a smooth projective variety faithfully flat over a quasi-compact quasi-separated scheme $S$ and suppose that there exists $s\in S$ such that $X_{s}$ has either an ample canonical bundle or ample anti canonical bundle. Let $(\boxtimes,\Ox_{X})$ be an $S$-linear tensor triangulated category structure on $Perf(X)$ which is geometrically reasonable and reasonable with respect to the adjunction pair $I_{s}^{\ast}\dashv I_{s}\,_{\ast}$. If $\omega_{X_{s}}$ is $\boxtimes_{I_{s}^{\ast}I_{s}\,_{\ast}}$-invertible, and $I_{s}^{\ast}(\mathbb{1})\simeq \Ox_{X_{s}}$, then there exists an open affine subset $U\subseteq S$ such that $\boxtimes_{U}$ coincides on objects with $\dteee$.
    \begin{proof}
    Let us write $\bGamma$ for the MTC structure enhancing $\boxtimes$. By Theorem \ref{thm:mainthm1} we have that $\boxtimes_{s}$ coincides on objects with $\dteee_{X_{s}}$. By Corollary \ref{thm:fiberisq} we have that there exists an open affine subset $U\subseteq S$ such that the 2-fold dg-bimodule $\bGamma_{U}$ enhancing $\boxtimes_{U}$ coincides with $\dteee_{X_{U}}$. 
    \end{proof}
\end{thm}
Our main corollary follows,
\begin{cor}\label{cor:maincor}
Suppose $X\to S$ is a smooth projective variety faithfully flat over a quasi-compact and quasi-separated scheme $S$, and assume that for every $s\in S$ the fiber $X_{s}$ has ample (anti-)canonical bundle. Furthermore suppose that $Spec(\boxtimes)$ is a smooth projective variety, faithfully flat over $S$ such that $\der[X]\simeq \der[Spec(\boxtimes)]$ then $X\cong Spec(\boxtimes)$ as S-schemes.
\begin{proof}
As $S$ is quasi-compact quasi-separated we can reduce the problem to a sufficiently small affine subset $U\subseteq S$ and suppose there is an equivalence $\der[Spec(\boxtimes_{U})]\simeq \der[X_{U}]$. By Theorem \ref{thm:mainthm22222222} we have a 2-fold dg-bimdoule $\bGamma_{U}$ enhancing $\boxtimes_{U}$ and which we know coincides with $\dteee_{X_{U}}$. Furthermore we can pick the equivalence $\der[X_{U}]\simeq \der[Spec(\boxtimes_{U})]$ to send $\omega_{X_{U}}$ to $\omega_{Spec(\boxtimes_{U})}$ and so this latter object is $\boxtimes_{U}$-invertible. This implies by \cite[Theorem 3.25]{toledo2024tensor} that $X_{U}\cong Spec(\boxtimes_{U})$. 
\end{proof}
\end{cor}
\section{Examples}\label{sec:examples}
We dedicate this section to sketch some possible interesting examples and applications of the constructions provided in this article. 
\begin{exmp}\label{example:fm}
    The first immediate example we ought to mention comes from the motivating situation that gives rise to the constructions described in this article. \\
    If $X\to S$ is a smooth projective variety faithfully flat over a quasi-compact quasi-separated base scheme $S$, then the derived category $\der[X]$ carries a tensor triangulated category structure $(\dteee_{X},\Ox_{X})$ which is linear over the usual sense and thus in the sense described in Section \ref{sec:prelim}. \\
    Additionally, for every non-trivial Fourier-Mukai partner $Y\to S$ of $X$, we obtain a non-equivalent tensor structure on $\der[X]$ by  means of transporting the structure $(\dteee_{Y},\Ox_{Y})$ on $\der[Y]$ to $\der[X]$ via the Fourier-Mukai transform equivalence. \\
    An interesting example which exhibits the need for the hypothesis in Theorem \ref{thm:mainthm22222222} comes from the classification due to Uehara in \cite{uehara2004example} where he shows:
    \begin{thm}
    \begin{enumerate}
        \item Let $p$ be a positive integer. Then there is a rational elliptic surface $S(p)$ such that $S(p)$ has a singular fiber of type $_{p}I_{0}$ and at least three non-multiple singular fibers of different Kodaira's types. 
        \item Let $N$ be a positive integer and $p$ a prime number such that $p>6(N-1)+1$. Then there are rational elliptic surfaces $T_{i}$, for $i\in\{1,\dots, N\}$ such that $T_{i}\not\cong T_{j}$ for $i\not= j$ and every $T_{i}$ is a Fourier-Mukai partner of $S(p)$. As a special case, $S=S(11)$ has an FM partner $T$ such that $T\not\cong S$. These $S$ and $T$ are birational, derived equivalent but not K-equivalent. 
    \end{enumerate}
    \end{thm}
    In particular for any of these surfaces $T_{i}$ we have that their general elliptic curve fibers $E$ all coincide. This implies that despite elliptic curves being too a class of varieties completely determined by their derived categories the conditions in Theorem \ref{thm:mainthm22222222} are necessary. In this case it is likely neither of the reasonableness conditions is satisfied, so that the transported structure doesn't send the unit $\Ox_{X}$ to an invertible object nor is it compatible with the pullback and pushforwards induced by (fiber products of) open subsets of $\mathbb{P}^{1}$ for any given fixed surface $T_{i}\to \mathbb{P}^{1}$.
\end{exmp}
\begin{exmp}\label{example:quiver}
    In the introduction of this work, we mentioned the results of \cite{liu2013recovering} to justify the investigation of tensor structures since their main results provides plenty of examples of tensor structures $(\der[X],\boxtimes,\mathbb{1})$ for varying smooth projective varieties $X$ which are not induced by tensor structures of the form $(\dteee_{Y},\Ox_{Y})$ for any smooth projective variety $Y$. \\
    Namely, they show that whenever there is an equivalence $\der[X]\simeq \der[\mr{Rep} (Q,R)]$ where $\mr{Rep} (Q,R)$ is the abelian category of finite representations of a quiver $Q$ with relations $R$ suitably compatible with the derived tensor product $\dteee_{(Q,R)}$ of quiver representations, then $Spec(\dteee_{(Q,R)})$ has an underlying topological space which is in bijection with the set of vertices of $Q$, and so one can conclude for example by Balmer's reconstruction theorem that the the tensor structure $\dteee_{(Q,R)}$ is not equivalent to $\dtee$. \\
    An equivalence of this nature exists for example as soon as the derived category $\der[X]$ admits a full strong exceptional collection $\{E_{i}\}$. Let us consider for example the Hirzebruch surface $\Sigma_{n}$ which can be obtained as the projective bundle associated to the rank 2 vector bundle $\Ox\oplus \Ox(-n)$ of $\mathbb{P}^{1}$. \\
    As such, it comes equipped with the projection structure map $\pi:\Sigma_{n}\to \mathbb{P}^{1}$. \\
    By Orlov's formula for projective bundles (\cite[Theorem 2.6]{orlov1992projective}) we can obtain a full exceptional collection, which can be furthermore shown to be strong:
    \begin{equation}\label{eq:excSn}
        \{\F,\F(1,0),\F(1,0),\F(1,1)\}.
    \end{equation}
    As before, let us consider the derived tensor product $\dteee_{(Q,R)}$ on $\der[\Sigma_{n}]$ obtained from the tensor product of representations of the quiver $Q$ induced by the exceptional collection \ref{eq:excSn}. \\
    If we denote by $U(m)$ the unitary representation defined by having the base field $k$ on the m-th vertex and $0$ on every other vertex and on every edge, one can see that any object $V\in \mr{Rep} (Q,R)$ in $Q$ is an extension of $U(m)$'s. \\
    If we let $E^{\bullet}\in \der[\mathbb{P}^{1}]$, one then can see that to calculate $\pi^{\ast}E^{\bullet}\dteee_{(Q,R)}\_$ one only needs to check at the representations $U(m)$. From this we obtain that $\dteee_{(Q,R)}$ is linear over $\dteee_{\mathbb{P}^{1}}$.
    In this sense, the tensor product $\dteee_{(Q,R)}$ again serves as an example of a monoidal structure not induced by a Fourier-Mukai transform. 
\end{exmp}
Our final remark on this direction comes from the following construction which can be found in Subotic's PhD thesis (\cite{subotic2010monoidal}). Let us briefly expand on this. We direct the reader interested in precise definitions and detailed constructions to Subotic's original work.
\begin{exmp}\label{ex:fukaya}
    The Homological Mirror Symmetry conjecture of Kontsevich establishes an equivalence between the derived category of coherent sheaves of a Calabi-Yau variety $X$ and the derived category of the Fukaya category of its mirror $X^{\vee}$. A large body of work has been dedicated to verifying this equivalence in particular instances. \\
    Let $M\to B$ be a Lagrangian torus fibration over a compact base $B$. Roughly speaking, the generalized Donaldson-Fukaya category $Don^{\#}(M)$ of $M$ is a category which has generalized Lagrangian correspondences $(L_{0},\dots,L_{m})$ with local systems $(S_{0},\dots,S_{m})$ as objects, and morphisms between two such objects are defined through the Floer homology with coefficients in the local system of their concatenation. See \cite[Definition 4.0.27]{subotic2010monoidal} for a precise definition. \\
    The main result in \cite{subotic2010monoidal} consists in the construction of a tensor structure $\hat{\otimes}$ on $Don^{\#}(M)$ which in the case $M$ is a torus, is shown to be compatible with the Homological Mirror Symmetry equivalence for elliptic curves of Polischuk-Zaslow in \cite{polishchuk1998categorical} when one equips the derived category of the associated elliptic curve $E$ with the derived tensor product $\dteee_{E}$. \\
    The surfaces from Example \ref{example:fm} are of non-zero Kodaira dimension and so we cannot use the existence of non-trivial Fourier-Mukai transforms in this setting to conclude Subotic's theorem does not apply to general mirror symmetry of toric fibrations. \\
    It would be of great interest to apply the analysis we have made of the sheaves of 2-fold dg-bimodules parametrizing bifunctors $\boxtimes:\der[X]\times \der[X]\to \der[X]$ and both how this parametrization behaves at fibers and neighborhoods of them as we can incorporate Subotic's monoidal equivalence in this framework. \\
    This point of view is to be thought as being in line with SYZ mirror symmetry which predicts the mirror behaviour in terms of toric fibrations. \\
    \end{exmp}
\addtocontents{toc}{\protect\setcounter{tocdepth}{0}}
\section*{Acknowledgements}
The first author is supported by a grant from Beijing Institute of Mathematical Sciences and Applications (BIMSA), and by the National Program of Overseas High Level Talent. The second author would like to thank Carlos Simpson for early conversations on some of the constructions described in the present article, and would also like to express his gratitude to the Instituto de Matemáticas Unidad Oaxaca of the National Autonomous University of Mexico for kindly providing him with office space where part of the first drafts of this work were written. 
\endgroup
\hypersetup{linkcolor=red}
\bibliographystyle{alpha}
\bibliography{RelativeMonoidalBO}
\end{document}